\documentclass[final]{siamltex}

\usepackage{latexsym, enumerate}
\usepackage{eepic}
\usepackage{epic}
\usepackage{graphicx}
\usepackage{color}
\usepackage{ifpdf}
\usepackage{amssymb,amsmath,epsfig, algorithm,algorithmicx,multirow}
\usepackage{amsfonts, dsfont}
\usepackage{subfigure}
\usepackage{multirow}
\usepackage{makecell}
\usepackage{bm}

\newcommand{\bb}{\mathbb}

\newcommand{\bt}{\bm{\theta}}
\newtheorem{thm}{Theorem}[section]
 \title{ Ensemble-based implicit sampling for Bayesian inverse problems with non-Gaussian priors
\thanks{L. Jiang acknowledges the support of Chinese NSF 11871387 .
}}

\author{Yuming Ba
\thanks{College of Mathematics and Econometrics, Hunan University, Changsha 410082, China. Email:yumingb@hnu.edu.cn.}
\and
Lijian Jiang\thanks{School of Mathematical Sciences, Tongji University, Shanghai 200092, China. Email: ljjiang@tongji.edu.cn. Corresponding author}
}
\begin{document}
%\title{ \large\bf Ensemble-based implicit sampling for non-Gaussian priors in Bayesian inverse problems}
%
%\author{
%Yuming Ba
%\thanks{College of Mathematics and Econometrics, Hunan University, Changsha 410082, China. Email:yumingb@hnu.edu.cn.}
%\and
%Lijian Jiang\thanks{School of Mathematical Sciences, Tongji University, Shanghai 200092, China. Email: ljjiang@tongji.edu.cn. Corresponding author}
%}

\maketitle

\begin{abstract}
In the paper, we develop  an ensemble-based implicit sampling method for Bayesian inverse problems.
For Bayesian inference,  the iterative ensemble smoother (IES)  and implicit sampling are integrated  to obtain importance ensemble samples, which build an importance density.
The proposed method shares a similar idea to importance sampling.  IES is used to  approximate  mean  and covariance of a posterior distribution. This provides the MAP point and the inverse of Hessian matrix,
which are necessary to construct the implicit map in implicit sampling. The importance  samples are generated by the implicit map and the corresponding weights are the ratio between the importance density and posterior density.
In the proposed method, we use the ensemble samples of IES to find the optimization solution of likelihood function and the  inverse of  Hessian matrix.  This approach avoids the   explicit computation for Jacobian matrix
and Hessian matrix, which are very computationally expensive in high dimension spaces.  To treat non-Gaussian models, discrete cosine transform  and Gaussian mixture model  are used to characterize
the non-Gaussian priors.  The ensemble-based implicit sampling method is extended to the non-Gaussian priors for exploring the posterior of unknowns in inverse problems.
The proposed method is  used for each individual  Gaussian model in the Gaussian mixture model. The proposed approach  substantially improves  the applicability of  implicit sampling method.
 A few numerical examples are presented to demonstrate  the efficacy of the proposed method with applications of inverse problems for subsurface flow problems  and anomalous diffusion models in porous media.

\end{abstract}



\pagestyle{myheadings}

\thispagestyle{plain}
\markboth{Y. Ba and L. Jiang}{Ensemble-based implicit sampling}

\section{Introduction}
The model inputs, such as parameters, the source and system structure, are often unknown in  practical models. For example, in porous media application, the practical models include  subsurface flows and anomalous diffusion models \cite{Wyss1986fractional}.  The unknown model inputs can be identified by integrating  noise observations and the prior information. The problem of identifying the unknown inputs in mathematical models have been intensively studied in the framework of inverse problems and many methods have been proposed \cite{Akella2011reservoir, Miller2013coefficient, chang2017surrogate}. The  inverse problem can be generally considered as an optimization problem, which minimizes the misfit between simulated observations and true observations.  The inverse problem  is often ill-posed. To avoid the ill-posedness, incorporating the penalty into the objective function is often necessary \cite{Hansen1992analysis}. In the paper, we use Bayesian inference to solve the inverse problems. Bayesian inversion  can provide not only the point estimate, but also the statistical distribution of unknowns and the prediction intervals of state variables.

In Bayesian inference, the posterior is often  concentrated in a small portion of the entire prior support. To accelerate the posterior exploration, the unknown inputs can be estimated by the samples from  a high probability region.   As a sampling method for Bayesian inference,  importance sampling  generates samples from a probability density function (pdf), which is only up to a multiplicative constant. The importance sampling method based on Monte Carlo has been widely used in inverse problems \cite{Mosegaard1995monte, Li2015adaptive}. It can be used as a method of sensitivity analysis, as an alternative of acceptance-rejection sampling and as the foundation to computing normalizing constants of probability densities. In dynamic systems, it is also an essential prerequisite for sequential Monte Carlo \cite{Doucet2000sequential}. However, the importance sampling is different from the standard Monte Carlo method,  which generates samples with  equal weights. The weights of importance sampling are from the proposal density and unequal. To obtain effective weights, the selection of the proposal density is critical. The proposal density is also called the importance density.
%For constructing an appropriate importance density, Sequential Monte Carlo is proposed for dynamic or nonlinear models in \cite{Doucet2000sequential, De2012}.

In the paper, we construct  the importance density through an  implicit sampling (IS) method \cite{chorin2009implicit}, which provides a data-informed importance density. The main idea of IS is to locate the high probability region and generate samples around the  Maximum A Posteriori (MAP) point. In  IS method,  it is required to compute  the MAP point of the posterior and the corresponding Hessian matrix for the negative logarithmic of the posterior.  There are many techniques to estimate  the MAP point, such as Markov chain Monte Carlo (MCMC) method \cite{Jiang2018bayesian, Efendiev2006preconditioning}, variational method \cite{Jin2010hierarchical}, and ensemble-based method \cite{Iglesias2013evaluation, ba2018tWO}. These methods are based on the Bayesian framework and provide a statistic analysis, which can give the prediction and credible intervals, pdf  and other statistical information. The inverse of Hessian matrix can be approximated by the posterior covariance matrix. Then the importance weights can be obtained by solving an implicit equation. The effectiveness of IS for Bayesian inverse problems has been studied  in \cite{chorin2010implicit, morzfeld2012random, Atkins2013implicit}.

The ensemble-based method, such as ensemble Kalman filter (EnKF) \cite{Evensen1994sequential, ba2018tWO} and ensemble smoother (ES) \cite{chang2017surrogate}, was proposed for data assimilation. In recent years,  the ensemble-based method has been used to forecast the state and estimate the unknown parameters. In the work, we  use an ensemble-based method to find the MAP point and approximate the inverse of Hessian matrix by ensemble samples. EnKF is widely used in Bayesian data assimilation  but brings  the problem of inconsistency. Although ES has no inconsistency issue, it  is a global update by assimilating all observations simultaneously and may perform poorly due to the single update. To improve ES method, Chen and Oliver proposed the iterative ensemble smoother (IES) \cite{chen2012ensemble}, which can deal with the high-dimensional and nonlinear problems in \cite{chen2013levenberg}. IES is still a Gaussian approximation, which can provide the first-order  and second-order moments of the posterior distribution. This brings  difficulty  for using IES to non-Gaussian distribution.

To treat  non-Gaussian models in the ensemble-based method, we can use two approaches:  parameterization  and non-parameterization. The parameterization approach uses a transform to gain the latent variables. The typical examples include  truncated pluri-Gaussian (TPG) \cite{Astrakova2015conditioning}, level set \cite{Mannseth2014relation}, multiple-point simulation (MPS) \cite{Strebelle2002conditional} and discrete cosine transform (DCT) \cite{JAFARPOUR2008HISTORY}.  The latent parameters can be easy to update by the ensemble-based method. In a practical situation,  we may only know the discretization of the physical domain for  a random field  with unknown  covariance information. For this situation,  DCT based on a Fourier-based transformation can parameterize the random field. It only depends on the discretization of the physical domain. DCT roots in the image processing and has been widely used for image compression. To construct DCT expansion, the cosine functions can be used to form a set of mutually orthogonal basis functions.
By overcome the challenge of the possible high dimensional parameters in DCT, we can use a truncated DCT, where the low frequency basis functions are retained and the high frequency are abandoned.
Truncated DCT can both capture the main features of the random  field  and improve the computation efficiency in Bayesian inversion.

Non-parameterization is another approach  for handling the non-Gaussian priors. One of the non-parameterization methods is semi-parametric, such as the mixture of  distributions, which adopts a convex combination of several distributions to approximate a posterior distribution. The common mixtures are  Beta \cite{Barron1999consistency}, triangular \cite{Perron2001bayesian} and Gaussian \cite{Roeder1997practical, li2016gaussian} mixture models. In the mixture model, each distribution can be seen a basis function. Thus, the goal is to find a set of optimal basis functions. The corresponding coefficients are the model weights, which imply the importance of distributions. In the paper, we use  Gaussian mixture model (GMM) for Bayesian inversion with non-Gaussian priors.  In particular, GMM can be coupled with EnKF to estimate the multimodality state distributions \cite{li2016gaussian}. When using GMM, the weight, mean and covariance of each Gaussian model become  unknown parameters, which need to be estimated. The expectation-maximization (EM) \cite{Chen2010demystified} method is available to forecast these GMM parameters. To obtain a good approximation, the number $k$ of models may be  large enough. However, large $k$ will bring the singularity of covariance, which results from the unbounded logarithmic form of the mixture. To avoid the singularity issue,  the Bayesian Ying Yang (BYY) harmony learning based on the general statistical learning framework has been proposed in \cite{Xu2001byy}. A model selection criteria and automatic model selection method, called BYY harmony data smoothing learning model selection criterion (BYY-HDS) \cite{hu2004investigation}, can be derived from BYY harmony learning to estimate GMM parameters. Here, the $k$ can be unknown and selected by BYY-HDS when the initial value of $k$ is large enough.

The goal of this paper is to combine IES method with IS to develop an ensemble-based implicit sampling for Bayesian inverse problems.  In the proposed method,  we use IES to compute   mean  and covariance.
The mean can be used as the MAP point and the covariance as the approximation of the inverse of Hessian matrix. Then an implicit map is  given by the mean and the Cholesky factorization of the covariance. The  implicit map
gives the importance ensemble samples, where the corresponding weights are the ratio between the importance density and posterior density. Resampling based on these weights may be necessary  to avoid the ensemble degeneration. Thus, IES provides the MAP point and the inverse of Hessian matrix, and then using IS generates the importance samples. For convenience, we refer the proposed method as IES-IS.
 In the paper, we apply IES-IS to non-Gaussian models based on  DCT and GMM, which are  used to handle the non-Gaussian priors.   When using BYY-HDS based GMM method, ensemble samples will be the training data to forecast the mean, covariance and weight of GMM.  IES-IS  is  performed for each Gaussian model in GMM.   This substantially improves  the applicability of the IES-IS method for Bayesian inference.
For some complex structures in the target  field, it may be not enough to capture the main features by the immediate estimation of the proposed method. To this end, we use  the post-processing based on the regularization \cite{Grebennikov2006fast}  to improve the connectivity of main features. In general, the penalty term can be given by a quadratic form, which can achieve the global minimum with the convex constraint in the closed interval.

The rest of the paper is organized  as follows. We begin with the general framework of implicit sampling for Bayesian inverse problems. In Section \ref{non-Gaussian}, we focus on the non-Gaussian priors, which can be handled using DCT and BYY-HDS based GMM. Section \ref{IES-IS-non} is devoted to developing  IES-IS based on DCT and GMM. In Section \ref{exam}, a few numerical examples are presented to illustrate the performance of the proposed method with applications of inverse problems for subsurface flow problems  and anomalous diffusion models in porous media. In particular, we recover  channel structures and  fractures in  porous media.  Some conclusions and comments are made finally.

%%%%%%%%%%%%%%%%%%%%%%%%%%%%%%%%%%%%%%%%%%%%%%%%%%%%%%%%%%%%%%%%%%%%%%%%%%%
%%%%%%%%%%%%%%%%%%%%%%%%%%%%%%%%%%%%%%%%%%%%%%%%%%%%%%%%%%%%%%%%%%%%%%%%%%%%%%%

\section{Ensemble-based IS for Bayesian inverse problems}
We assume that a model problem is defined as
\[
\mathcal{U}(u;a(x))=f(x),
\]
where $\mathcal{U}$ is a generic forward operator and describes the relation of the coefficient $a(x)$, state $u$ and source term $f(x)$. For the Bayesian inverse problem, $a(x)$ and $f(x)$ may be unknown and assume to be characterized by  $a(x,\bt)$ and  $f(x,\bt)$, respectively,  in a finite dimensional parameter space, where $\bt\in\bb R^{N_{\theta}}$ is the unknown parameter. Let $g$ be the forward operator mapping the model parameter $\bt$ to the observation space, i.e., $y=:g(\bt)\in\bb R^{N_d}$.
Then the observation model can be given by
\[
\bm d=g(\bt)+\bm\varepsilon,
\]
where $\bm\varepsilon\in\bb R^{N_d}$ is the observation noise. In the paper, we assume that $\bm\varepsilon$ is independent of $\bt$ and $\varepsilon\sim N(\bm 0,\bm C_D)$. From the observation model, we have the likelihood function
\begin{equation}\label{Lik}
p(\bm d|\bt)\propto\exp\bigg[-\frac{1}{2}\bigg({\bm d}-g({\bt})\bigg)^T{\bm C_D}^{-1}\bigg({\bm d}-g({\bt})\bigg)\bigg].
\end{equation}

Given a prior $p(\bt)$, the conditional posterior density function can be derived by  Bayes rule
\[
p(\bt|\bm d)=\frac{p(\bt)p(\bm d|\bt)}{\int p(\bt)p(\bm d|\bt)d\bt},
\]
where $\int p(\bt)p(\bm d|\bt)d\bt$ is a constant independent of $\bt$.
Let $\bb{F}({\bt})=-\log \bigg(p({\bt})p({\bm d}|\bt)\bigg)$.
The goal of Bayesian inverse problem is to find a solution to minimizing  $\bb F(\bt)$, i.e.,
\begin{equation}\label{Bcost}
\widehat{\bt}=\arg\min_{\bt\in{\bb R}^{N_\theta}} \bb F(\bt).
\end{equation}
Thus,  $\widehat\bt$ is the MAP point of $p(\bt|\bm d)$.
In the Bayesian framework, $\widehat\bt$ can be approximated by the expectation of $\bt$ with respect to $p(\bt|\bm d)$
\begin{equation}\label{IS}
\widehat\bt\approx{\bb E}[\bt]=\int \bt p(\bt|\bm d)d\bt,
\end{equation}
where $\bb E[.]$ is the expectation operator.

\subsection{Bayesian inference using importance sampling}
We can use  Monte Carlo method by drawing $N_e$ independent samples from $p({\bt}|\bm d)$ to approximate the integral in equation (\ref{IS}).  Monte Carlo integration is often used  if we can sample from the target distribution. However, drawing samples from the target distribution is often difficult in practice.  Thus, we want  to seek an alternative distribution, which can be easy to sample. This motivates the importance sampling,  where the idea is to draw samples from a proposal distribution and re-weight the integral using the importance weights such that  a proper  distribution is targeted. The proposal density is also called the importance density. The importance sampling can bring enormous gains, making an otherwise infeasible  problem amenable to Monte Carlo.

When drawing samples from $p(\bt|\bm d)$ is infeasible, we need to find a importance density function $q(\bt|\bm d)$ to replace it.  Then we  have
\[
{\bb E}[\bt]=\int \frac{\bt p(\bt|\bm d)}{q(\bt|\bm d)}q(\bt|\bm d)d{\bt}={\bb E}_q[w({\bt})\bt],
\]
where $w({\bt})=\frac{p({\bt}|\bm d)}{q({\bt}|\bm d)}$ and ${\bb E}_q[\cdot]$ denotes the expectation with respect to $q({\bt}|\bm d)$. Then we sample  from $q$ instead of $p$.  Here  the adjustment factor $w$ is called the likelihood ratio. Thus the importance sampling estimate of (\ref{IS}) is given by
\[
\widehat{{\bb E}}_q[w({{\bt}})\bt]=\frac{1}{N_e}\sum_{i=1}^{N_e}w({\bt}_i)\bt_i,
\]
where ${\bt}_i$ is drawn from $q$. In Bayesian inference, $p({\bt}|\bm d)$ is only up to a normalizing constant, i.e., $w({\bt})=cw^0({\bt})$, where $w^0({\bt})$ can be obtained but $c$ is unknown. For this case, we compute the ratio estimate
\begin{equation*}
\widetilde{{\bb E}}_q[w({\bt})\bt]=\frac{\sum_{i=1}^{N_e}w({\bt}_i)\bt_i}{\sum_{i=1}^{N_e}w({\bt}_i)}
\end{equation*}
instead of $\widehat{{\bb E}}_q[w({\bt})\bt]$. Thus constructing the effective importance density function $q$ is critical  in Bayesian inverse problem.

%%%%%%%%%%%%%%%%%%%%%%%%%%%%%%%%%%%%%%%%%%%%%%%%%%%
\subsection{Implicit sampling method}
Implicit sampling \cite{chorin2009implicit} generates samples  by  an implicit map. It can construct the importance probability density function and ensure the efficacy of importance sampling. To implement IS, we first need to compute the MAP point of $p(\bt|\bm d)$ and Hessian matrix of $\bb F(\bt)$, and then generate samples from the high probability region of the posterior density.
Assume that  the minimum of $\bb F(\bt)$  exists. Let
\[
\varphi_{\bb F}=\min {\bb F}\quad  \text{and} \quad  \widetilde{{\bm\mu}}=\arg\min {\bb F}.
\]
We first find the high probability region of posterior density function by minimizing $\bb F({\bt})$. Then we generate samples around $\widetilde{{\bm\mu}}$.

If $g(\cdot)$ is  nonlinear, the posterior density function $p(\bt|\bm d)$ may be non-Gaussian, even though the prior $p(\bt)$ is Gaussian. Thus sampling from $p(\bt|\bm d)$ may be difficulty. In IS,
we choose  a reference random variable $\bm\xi$ with probability density function $b({\bm\xi})\propto e^{-\bb B({\bm\xi})}$, which is easy to sample (the reference random variable is  Gaussian in the paper). We assume that
the minimum of $\bb B(\bt)$ exists and $\varphi_{\bb B}=\min{\bb B}$. To generate the samples of ${\bt}$, we draw samples from $b({\bm\xi})$ and then solve the following implicit equation
\begin{equation}\label{imequation}
\bb F({\bt})-\varphi_{\bb F}=\bb B({\bm\xi})-\varphi_{\bb B}.
\end{equation}
The implicit equation can be solved by many approaches, such as  random map \cite{morzfeld2012random},  linear map \cite{chorin2010implicit} and the connection with optimal map \cite{el2012bayesian}. For the different  methods, the resulting samples may have different weights.

In the paper, we focus on  the linear map, where $\bb F({\bt})$ is approximated by
\begin{equation}\label{limap}
\bb F({\bt})\approx{\bb F}_0({\bt}):=\varphi_{\bb F}+\frac{1}{2}({\bt}-\widetilde{{\bm\mu}})^T{\bm H}({\bt}-\widetilde{{\bm\mu}}),
\end{equation}
where $\bm H$ is the Hessian matrix at $\widetilde{{\bm\mu}}$. If  we take $\bm\xi \sim N(\bm 0,\bm I)$, where $\bm I$ is the identity matrix. Then $\varphi_{\bb B}=0$ and $\bb B(\bm\xi)=\frac{1}{2}\bm\xi^T\bm\xi$.  By equation (\ref{imequation}) and (\ref{limap}), we get the approximated implicit equation
\begin{equation}\label{aime}
\bb F_0({\bt})-\varphi_{\bb F}=\frac{1}{2}{{\bm\xi}}^T{\bm\xi}.
\end{equation}
Let  $\bm L$ be the Cholesky factorization of $\bm H^{-1}$.
Then
\[
{\bt}=\widetilde{{\bm\mu}}+{\bm L}\bm\xi
\]
solves the equation  (\ref{aime}).
The corresponding weight of sample $\bt$ is
\[
w\propto \exp(\bb F_0({\bt})-\bb F({\bt})).
\]

%%%%%%%%%%%%%%%%%%%%%%%%%%%%%%%%%%%%%%%%%%%%%

\subsection{IS based on the iterative ensemble smoother}\label{IES-ISpr}
Iterative ensemble smoother (IES) was proposed in \cite{chen2012ensemble}, which is an ensemble method for Bayesian inverse problems. One of the IES methods is the modified Levenberg-Marquart method for  ensemble randomized maximum likelihood (LM-EnRML) \cite{chen2013levenberg}. For LM-EnRML method, a modification is made to approximate the inverse of Hessian matrix such that the explicit computation of the Jacobian matrix of $g$ is avoided.
We couple IES with IS to develop the ensemble-based implicit sampling method, where  IES is used  to obtain the MAP point of $p(\bt|\bm d)$ and the inverse of Hessian matrix of $\bb F(\bt)$, and IS is used to  generate  high probability samples.

From equation (\ref{Lik}), the likelihood function belongs to the exponential family. Thus we use a mixture of natural conjugate priors to approximate any prior $p(\bt)$ \cite{Dalal1983approximating}. In the paper, we use Gaussian density functions as the natural conjugate prior. Then we have
\[
p(\bt)=\sum_{i=1}^k\pi_i p(\bt|\bt^{pr,i},\bm C_{\theta,i}),\quad \pi_i>0\quad \text{and}\quad\sum_{i=1}^k\pi_i=1,
\]
where $p(\bt|\bt^{pr,i},\bm C_{\theta,i})$ is a Gaussian density function with the mean $\bt^{pr,i}$ and covariance $\bm C_{\theta,i}$,  and $\pi_i$ is the mixing probability. Thus the posterior density function
\[
p(\bt|\bm d)\propto\bigg(\sum_{i=1}^k\pi_i p(\bt|\bt^{pr,i},\bm C_{\theta,i})\bigg)p(\bm d|\bt).
\]
Let $p_i(\bt|\bm d)=\frac{p(\bt|\bt^{pr,i},\bm C_{\theta,i})p(\bm d|\bt)}{\int p(\bt|\bt^{pr,i},\bm C_{\theta,i})(\bt)p(\bm d|\bt)d\bt}$, i.e., $p_i(\bt|\bm d)\propto p(\bt|\bt^{pr,i},\bm C_{\theta,i})p(\bm d|\bt)$.
By  the convex combination of posterior density function,
\[
\bb E_p[\bt]=\sum_{i=1}^k{\pi}_i\bb E_{p_i}[\bt].
\]
%\begin{eqnarray*}
% % \nonumber to remove numbering (before each equation)
%\bb E_p[\bt]&=&\int\frac{\bt(\sum_{i=1}^k\pi_i l_i(\bt)p(\bm d|\bt))}{ \sum_{i=1}^k\pi_i c_i}d\bt\\
%&=&\sum_{i=1}^k\frac{\pi_i c_i\bb E_{p_i}[\bt]}{\sum_{i=1}^k\pi_i c_i}=\sum_{i=1}^k\tilde{\pi}_i\bb E_{p_i}[\bt],
%\end{eqnarray*}
%where $\tilde\pi_i=\frac{\pi_i c_i}{\sum_{i=1}^k\pi_i c_i}$.
To compute the expectation of ${\bb E}_{p_i}[\cdot]$, we need to find the posterior density function $p_i(\bt|\bm d)$. For the ensemble method, $p_i(\bt|\bm d)$ can be constructed by ensemble samples.

Let
\begin{eqnarray*}
\bb F_i(\bt)&=&-\log \bigg(p(\bt|\bt^{pr,i},\bm C_{\theta,i})p(\bm d|\bt)\bigg)\\
&=&\frac{1}{2}\bigg({{\bm d}}-g({\bt})\bigg)^T{\bm C}_D^{-1}\bigg({{\bm d}}-g({\bt})\bigg)+\frac{1}{2}({\bt}-{\bt}^{pr,i})^T{\bm C}_{\theta,i}^{-1}({\bt}-{\bt}^{pr,i})+c,
\end{eqnarray*}
 where $g(\cdot)$ the forward operator and ${\bm C}_D$ is the covariance matrix of observation error, and
$c$ is a constant independent of $\bt$. The ensemble samples can be considered as the minimizer of $\bb F_i(\bt)$. We note that  minimizing $\bb F_i(\bt)$ is equivalent to minimizing  the following function
\begin{equation}\label{cost}
\bb W_i(\bt)=\frac{1}{2}\bigg({{\bm d}}-g({\bt})\bigg)^T{\bm C}_D^{-1}\bigg({{\bm d}}-g({\bt})\bigg)+\frac{1}{2}({\bt}-{\bt}^{pr,i})^T{\bm C}_{\theta,i}^{-1}({\bt}-{\bt}^{pr,i}),\quad i=1,\cdots,k.
\end{equation}
 For most practical applications, $g(\cdot)$ is nonlinear.
Minimizing equation (\ref{cost}) with all observations simultaneously is called ensemble smoother, which is just one-step iteration and inaccurate for the high-dimensional or nonlinear problems. Thus we devote to using the iterative scheme, which is called IES method.

Let  ${\bm G}_l$ be the Jacobian matrix at $\bt_l$ at $l$-th iteration step of model $i$.   For  Gauss-Newton method, the gradient of equation (\ref{cost}) can be expressed as
\[
\quad \nabla_{\theta}\bb W_i({\bt}_l)\approx\bigg[{\bm C}_{\theta,i}^{-1}({\bt}_l-{\bt}^{pr,i})+{\bm G}_l^T{\bm C}_D^{-1}\bigg(g({\bt}_l)-{\bm d}\bigg)\bigg]
\]
 and Hessian matrix $\bm H_i\approx{\bm C}_{\theta,i}^{-1}+{\bm G}_l^T{\bm C}_D^{-1}{\bm G}_l$.
To avoid the influence of large data mismatch in early iterations and accelerate the convergence, we modify the Hessian matrix by the Levenberg-Marquart method \cite{chen2013levenberg}.
Besides, two further modifications in \cite{le2016adaptive} are necessary to implement the iterative update formula. Let $\overline{\bm G}_l$ be  the Jacobian matrix at the ensemble mean $\bar{\bt}_l$, where ${\bar\bt}_l=\frac{1}{N_e}\sum_{j=1}^{N_e}\bt_{l}^j$ and $N_e$ is the ensemble size. We use ${\bm C}_{\theta_l}$ to replace ${\bm C}_{\theta,i}$ in Hessian matrix and $\overline{\bm G}_l$ to replace ${\bm G}_l$. For IS based on IES method, we have the iterative scheme
\begin{eqnarray*}
\widetilde{\bt}_{l+1}^j=&{\bt}_l^j-\bigg((1+\lambda){\bm C}_{\theta_l}^{-1}+\overline{\bm G}_l^T{\bm C}_D^{-1}\overline{\bm G}_l\bigg)^{-1}\bigg[{\bm C}_{\theta,i}^{-1}({\bt}_l^j-{\bt}^{pr,i})+\overline{\bm G}_l^T{\bm C}_D^{-1}\bigg(g({\bt}_l^j)-{\bm d}\bigg)\bigg]\\
=&{\bt}_l^j-\frac{1}{1+\lambda}\bigg[{\bm C}_{\theta_l}-{\bm C}_{\theta_l}\overline{\bm G}_l^T\bigg((1+\lambda){\bm C}_D+\overline{\bm G}_l{\bm C}_{\theta_l}\overline{\bm G}_l^T\bigg)^{-1}\overline{\bm G}_l{\bm C}_{\theta_l}\bigg]{\bm C}_{\theta,i}^{-1}({\bt}_l^j-{\bt}^{pr,i})\\
&-{\bm C}_{\theta_l}\overline{\bm G}_l^T\bigg((1+\lambda){\bm C}_D+\overline{\bm G}_l{\bm C}_{\theta_l}\overline{\bm G}_l^T\bigg)^{-1}\bigg(g({\bt}_l^j)-{\bm d}\bigg).
\end{eqnarray*}
Then the Kalman gain
\begin{equation}\label{IES_kal}
{\bm K}_l={\bm C}_{\theta_l}\overline{\bm G}_l^T\bigg((1+\lambda){\bm C}_D+\overline{\bm G}_l{\bm C}_{\theta_l}\overline{\bm G}_l^T\bigg)^{-1}.
\end{equation}
%In equation $(\ref{IES_kal})$, ${\bm C}_{\theta_l}\overline{\bm G}_l^T$ and
%$\overline{\bm G}_l{\bm C}_{\theta_l}\overline{\bm G}_l^T$ have different approximation forms for $k=1$ and $k\geq2$. We will give the concrete expression in Algorithm \ref{IES-IS-DCT} and \ref{IES-IS-GMM}.
The intermediate ensemble samples are generated by
\begin{equation}\label{IES}
\widetilde{\bt}_{l+1}^j={\bt}_l^j-\frac{1}{1+\lambda}({\bm C}_{\theta_l}-{\bm K}_l\overline{\bm G}_l{\bm C}_{\theta_l}){\bm C}_{\theta,i}^{-1}({\bt}_l^j-{\bt}^{pr,i})-{\bm K}_l\bigg(g({\bt}_l^j)-{\bm d}\bigg),\quad j=1,\cdots,N_e.
\end{equation}

For IS method, the inverse of modified Hessian matrix
\begin{equation}\label{IS_H}
\widetilde{{\bm H}}_{i}^{-1}\approx\frac{1}{1+\lambda}({\bm C}_{\theta_l}-{\bm K}_l\overline{\bm G}_l{\bm C}_{\theta_l}),
\end{equation}
and the MAP point
\begin{equation}\label{IS_mu}
{\widetilde{\bm\mu}}_{i}\approx\frac{1}{N_e}\sum_{j=1}^{N_e}\widetilde{\bt}_{l+1}^j.
\end{equation}
Assume that  the minimum of $\bb F_i(\bt)$ exists and $\phi_{\bb F_i}=\min\bb F_i(\bt)=\min\bb W_i(\bt)+c$. Let
$\phi_{\bb W_i}=\min\bb W_i(\bt)$ and $\bb W_0^i(\bt)=\phi_{\bb W_i}+\frac{1}{2}\bm\xi^T\bm\xi$. We apply  equation (\ref{limap}) and (\ref{aime}) to $\bb F_i(\bt)$ and have
%\[
%\bb F_i(\bt)=\bb W(\bt)+c\approx\bb F_0^i(\bt)=\phi_{\bb F_i}+\frac{1}{2}\bm\xi^T\bm\xi=\phi_{\bb W_i}+c+\frac{1}{2}\bm\xi^T\bm\xi.
%\]
%Let . We can get
\begin{equation}\label{iaime}
\bb W_i(\bt)\approx\bb W_0^i(\bt)=\phi_{\bb W_i}+\frac{1}{2}\bm\xi^T\bm\xi.
\end{equation}
The importance ensemble samples can be obtained by solving equation (\ref{iaime}), i.e.,
\begin{equation}\label{sam_IS}
\bt_{l+1}^{j,i}=\widetilde{{\bm\mu}}_{i}+{\bm L}_i\bm\xi_j
\end{equation}
with the weights
\[
w_j^i\propto \exp(\bb W_0^i(\bt_{l+1}^{j,i})-\bb W_i(\bt_{l+1}^{j,i})),\quad j=1,\cdots,N_e,\quad i=1,\cdots,k,
\]
where $\bm L_i$ is the Cholesky factorization of $\widetilde{\bm H}_i^{-1}$. For the discrete structure, the weights may approach  infinity because of the exponential growth. To ensure a proper size of the effective samples, we can scale the difference
$\bb W_0^i(\bt_{l+1}^{j,i})-\bb W_i(\bt_{l+1}^{j,i})$, where the order of weights retains unchanged. This does  not affect the solution of the implicit equation. Thus the importance samples can be generated by equation (\ref{sam_IS}) with the modified  weights
\begin{equation}\label{weight}
w_j^i\propto \exp(\frac{\bb W_0^i(\bt_{l+1}^{j,i})-\bb W_i(\bt_{l+1}^{j,i})}{\rho}),\quad j=1,\cdots,N_e,\quad i=1,\cdots,k.
\end{equation}
The selection of the scale parameter $\rho$ is carefully discussed in  \cite{sjz18}.  
To avoid ensemble collapse, we do  resampling and  redistribute the weights. For each model $\bb F_i(\bt)$, we  obtain  the importance ensemble samples $\{\bt_{l+1}^{j,i}\}_{j=1}^{N_e}$. In  the end, we combine $k$ ensembles using a membership probability matrix to get the update ensemble $\{\bt_{l+1}^j\}_{j=1}^{N_e}$.

\section{Priors based on  DCT and GMM}\label{non-Gaussian}
Ensemble-based method may not work well for the problems with non-Gaussian priors. In order to overcome the difficulty, we use suitable parameterization methods to characterize the non-Gaussian field.
In the paper, we focus on DCT and GMM to treat the non-Gaussian priors.  We apply the proposed IES-IS method to the priors described by  DCT and GMM.

\subsection{Parameterization based on DCT}\label{Pa}
The goal of parameterization methods is to obtain the latent variables by a transform, which can be updated by the ensemble-based method. The widely used parameterization methods are truncated pluri-Gaussian (TPG), level set, multiple-point simulation (MPS) and discrete cosine transform (DCT). To get the unknown parameter $\bt$, these parameterizations are used to the unknown input. In the paper, we focus on DCT  \cite{JAFARPOUR2008HISTORY}, which  is a Fourier-based transformation and can extract the important features of a random field in the Bayesian inverse problem. As known, paramerization by KLE  needs the mean and covariance information of the random field. But for  DCT method, the physical domain  discretization of the random field is enough to construct the basis functions. Besides, the separability of DCT basis makes the efficient computation of basis functions.

 Without loss of generality, we consider the unknown  function $a(x)$ defined in a two dimensional spatial domain, and $N_x\times N_y$ uniform grid is used to discretize the function.  Then $a(x)$  can be expressed as ${\bm A}(m,n)$ $(m=0,\cdots,N_x-1;n=0,\cdots,N_y-1)$. Thus the general forward DCT of the input field ${\bm A}(m,n)$ has the form
\begin{equation*}
  {\bt}_{N_m}(i,j)=\frac{2\alpha(i)\alpha(j)}{\sqrt{N_xN_y}}\sum_{m=0}^{N_x-1}\sum_{n=0}^{N_y-1}{\bm A}(m,n)\cos\bigg[\frac{\pi(2m+1)i}{2N_x}\bigg]\cos\bigg[\frac{\pi(2n+1)j}{2N_y}\bigg],
\end{equation*}
where $i=0,\cdots,N_x-1$, $j=0,\cdots,N_y-1$ and
\[
  \alpha(i)=
   \begin{cases}
    \frac{1}{\sqrt{2}}& i=0,\\
     1 & \text{otherwise}.
  \end{cases}
\]
In the inversion problem, we get $\bm A$ through estimating the unknown $\bt_{N_m}$. Thus the inverse 2D DCT is necessary for the inverse problem, as shown below
\begin{equation}\label{inDCT}
  {\bm A}(m,n)=\frac{2}{\sqrt{N_xN_y}}\sum_{i=0}^{N_x-1}\sum_{j=0}^{N_y-1}\alpha(i)\alpha(j){\bt_{N_m}}(i,j)\cos{\bigg[}\frac{\pi(2m+1)i}{2N_x}{\bigg]}\cos{\bigg[}\frac{\pi(2n+1)j}{2N_y}{\bigg]},
\end{equation}
where $m=0,\cdots,N_x-1$, $n=0,\cdots,N_y-1$.

Due to the separability property of DCT basis functions, (\ref{inDCT}) can be written as
\begin{equation}\label{seDCT}
  {\bm A}(m,n)=\sum_{i=0}^{N_x-1}\sqrt{\frac{2}{N_x}}\alpha(i)\bigg\{\sum_{j=0}^{N_y-1}\sqrt{\frac{2}{N_y}}\alpha(j){\bt_{N_m}}(i,j)\cos\bigg[\frac{\pi(2n+1)j}{2N_y}\bigg]\bigg\}\cos\bigg[\frac{\pi(2m+1)i}{2N_x}\bigg].
\end{equation}
We can implement DCT by a vector form, i.e., $\bm A$ and $\bt_{N_m}$ can be represented as the vectors. Thus
${\bm A}\in\bb R^{N_m}$ and $\bt_{N_m}\in\bb R^{N_m}$.
Let $\bm\Phi_{N_m} \in\bb R^{N_m\times N_m}$ denote basis function matrix with respect to $\bt_{N_m}$ and $\bb I_{N_m}$ denote a set of natural numbers given by
\[
\bb I_{N_m}=\{i*N_y+j\},\quad i\in\{0,\cdots,N_x-1\},\quad j\in\{0,\cdots,N_y-1\},
\]
where all the possible combinations for $i$ and $j$ are taken into account. The $r$th column of $\bm\Phi_{N_m}$ can be expressed as
\begin{equation}\label{DCT}
{\bm\phi}_r=\frac{2\alpha(i)\alpha(j)}{\sqrt{N_m}}
\left[
 \begin{array}{c}     % 2ÁУ¬¾ÓÖÐ
 \cos\bigg[\frac{\pi(2\times0+1)i}{2N_x}\bigg]\cos\bigg[\frac{\pi(2\times0+1)j}{2N_y}\bigg] \\
  \vdots \\
 \cos\bigg[\frac{\pi(2\times(N_x-1)+1)i}{2N_x}\bigg]\cos\bigg[\frac{\pi(2\times(N_y-1)+1)j}{2N_y}\bigg] \\
 \end{array} \right]\in\bb R^{N_m}.
\end{equation}
Then (\ref{seDCT}) becomes
\[
{\bm A}=\bm\Phi_{N_m}{\bt}_{N_m}=[{\bm\phi}_0,{\bm\phi}_1,\cdots,{\bm\phi}_{N_y-1},\cdots,{\bm\phi}_{N_m-1}]{\bt}_{N_m}=\sum_{r=0}^{N_m-1}{{\theta}}_r{\bm\phi}_r,
\]
where $N_m=N_xN_y$ and the subscript $r\in \bb I_{N_m}$ corresponding to a pair of indices $(i,j)$.

From equation (\ref{DCT}), we note that $\bm\Phi_{N_m}$ can be pre-computed and data-independent. Thus DCT basis functions only need to be calculated and stored once.
We choose  the low frequency basis functions in $\bm\Phi_{N_m}$ to retain the main features of a random field. Then we reduce the dimension of unknown parameters without losing the main  features. The selection of basis functions makes the low frequency basis functions retained and discards the high frequency basis funtions.
The truncated DCT expansion by the first $N_c$ terms can be represented by
\begin{equation}\label{TDCT}
{\bm A}=\sum_{r=0}^{N_c-1}{{\theta}}_r{\bm\phi}_r={\bm\Phi}{{\bt}},
\end{equation}
where $\bm\Phi\in\bb R^{N_m\times N_c}$ is the first $N_c$ columns of $\bm\Phi_{N_m}$ and $\bt\in\bb R^{N_c}$.
Using equation (\ref{TDCT}), the IES-IS  can be performed in the low-dimensional stochastic subspace.

\subsection{Nonparametric method based on GMM}
 Compared with the parameterizations, semi-parametric is a different approach. The typical method is to adopt the mixture models to approximate the unknown distribution. We assume that a mixture of distributions can be described as any convex combination of other distributions $P_i$,
\[
\sum_{i=1}^k\pi_iP_i({\bt}),
\]
where $\sum_{i=1}^k\pi_i=1$ and $k>1$, and $\{P_i\}_{i=1}^k$ are from a parametric family. For an unknown distribution, there is a trade-off between the perfect representation of the unknown distribution and the useful estimation of the mixture. To obtain a good approximation of the distribution, $k$ may be large enough. The mixture models can be considered  as using a few  basis distributions  to approximate the unknown distribution.

Many mixtures have been applied to the Bayesian inverse problems, such as Beta, triangular and Gaussian mixture models. In the paper, we focus on  Gaussian mixture model (GMM). We assume $\mathbb{L}({\bt})$ be a Gaussian mixture density function consist of a convex combination of $k$ Gaussian density functions, i.e.,
\[
\bb L({\bt})=\sum_{i=1}^k\pi_ip({\bt}|{\bm\mu}_i,{\bm\Sigma}_i).
\]
Here $p(.|{\bm\mu_i},{\bm\Sigma_i})$ is a Gaussian density function with the mean ${\bm\mu_i}$ and covariance ${\bm\Sigma_i}$. GMM parameters $\{\pi_i,{\bm\mu}_i,{\bm\Sigma}_i\}_{i=1}^k$ are unknown, which need to be identified in the Bayesian inverse problems. Let ${\bm q}=\{\pi_i,{\bm\mu}_i,{\bm\Sigma}_i\}_{i=1}^k$. When the estimation of $\bm q$ obtained, we have the probability density function with respect to $\bt$. The expectation of ${\bt}$ can be computed using the convexity
\[
\mathbb{E}[{\bt}]=\sum_{i=1}^k\pi_i{\bm\mu}_i.
\]

To get an accurate  estimation, we can update $\bm q$ by an iteration process. The widely used method is the Expectation-maximization (EM) for estimating GMM parameters with known $k$. For sufficiently large $k$, the covariance matrices may be singular in EM algorithm. This is an inherent problem that the logarithmic form of $\bb L({\bt})$ is unbounded. For estimating GMM parameters with unknown $k$, the Bayesian Ying Yang harmony data smoothing (BYY-HDS) learning model selection criterion is proposed in \cite{hu2004investigation}. To get a proper approximation, BYY-HDS method can automatically screen models by minimizing the function
\begin{equation}\label{BYY}
\bb J_{BYY-HDS}({\bm q}_k^h,k)=\sum_{i=1}^k\pi_i(\frac{1}{2}\log|{\bm\Sigma}_i|+\frac{1}{2}h^2{\text {tr}}[{\bm\Sigma}_i^{-1}]-\log\pi_i),
\end{equation}
where ${\bm q}_k^h=({\bm q}_k,h)$ and ${\text {tr}}[\cdot]$ denotes the trace operator of the matrix. Due to $k$ unfixed, $\bm q_k$ denotes the parameters of $k$ models at the current iteration. For convenience, EM method using BYY-HDS learning model selection criterion is called the smoothed EM (SmEM).

In the paper, we devote to using SmEM method to ensemble samples instead of the observation data, which is proposed in \cite{li2016gaussian}. For SmEM algorithm, we first need to set the maximum and minimum of $k$, where the maximum is large enough to automatically screen the mixture model and minimum is larger than 1.
At the $l$-th iteration step, we select and discard the models corresponding to the smaller weights $\pi_i$. Thus $k$ may decrease with respect to  the iterations. We update GMM parameters, which can be considered as E-step and M-step. E-step of SmEM can be expressed as
\begin{equation}\label{Estep}
\gamma_{i,j}=\frac{\pi_ip({\bt}^j|{\bm\mu}_i,{\bm\Sigma}_i)}{\sum_{i=1}^k\pi_ip({\bt}^j|{\bm\mu}_i,{\bm\Sigma}_i)},\quad i=1,\cdots,k,\quad j=1,\cdots,N_e.
\end{equation}
When we get the samples probability matrix $(\gamma_{i,j})$, M-step in SmEM can be represented by
\begin{equation}\label{Mstep}
\left\{
 \begin{aligned}
\pi_i&=\frac{1}{N_e}\sum_{j=1}^{N_e}\gamma_{i,j},\\
{\bm\mu}_i&=\frac{1}{N_e\pi_i}\sum_{j=1}^{N_e}\gamma_{i,j}{\bt}^j,\\
\bm\Sigma_i&=\frac{1}{N_e\pi_i}\sum_{j=1}^{N_e}\gamma_{i,j}({\bt}^j-{\bm\mu}_i)({\bt}^j-{\bm\mu}_i)^{T}+h_l^2{\bm I},\quad i=1,\cdots,k,\quad j=1,\cdots,N_e,\\
 \end{aligned}
 \right.
\end{equation}
where $h_l$ is the smoothing parameter at the $l$-th  iteration step. By equation (\ref{Mstep}), we note that ${\bm\Sigma}_i$ is modified to avoid the singularity, which often occurs in EM method. The  $h_l$ is critical for SmEM. We use an iteration scheme to update the smoothing parameter $h_l$, i.e.,
\begin{equation}\label{smooth}
h_{l+1}=h_{l}+\eta \bb S(h_{l}),\\
\end{equation}
where $\eta$ is the step length constant and
\[
\bb S(h_{l})=\frac{N_{\theta}}{h_{l}}-h_{l}\sum_{i=1}^k\pi_i{\text {tr}}[{\bm\Sigma}_i^{-1}]-\frac{\sum_{i=1}^{N_e}\sum_{j=1}^{N_e}\beta_{i,j}\|{\bt}^i-{\bt}^j\|^2}{h_{l}^3}\\
\]
with
\[
\beta_{i,j}=\frac{\exp(-0.5\frac{\|{\bt}^i-{\bt}^j\|^2}{h_{l}^2})}{\sum_{i=1}^{N_e}\sum_{j=1}^{N_e}\exp(-0.5\frac{\|{\bt}^i-{\bt}^j\|^2}{h_{l}^2})}.\\
\]
The detailed procedure  is presented in Algorithm \ref{SmEM}.
\begin{algorithm}
\caption{SmEM algorithm for estimating GMM parameters}
  \textbf{Input}: Given a set of data $\{{\bt}^j\}_{j=1}^{N_e}$, the initial ${\bm q}_k=\{{\pi}_i,{\bm\mu}_i,{\bm\Sigma}_i\}_{i=1}^k$,  sufficiently small\\
  $~~~~~~~~~~$ $\epsilon$, the initial $h_0=\sqrt{\frac{1}{N_{\theta}N_e^3}\sum_{i=1}^{N_e}\sum_{j=1}^{N_e}\|\bt^i-\bt^j\|^2}$, positive integer $k_{max}$ and\\
  $~~~~~~~~~~$  $k_{min}$ and given $\eta$.\\
  \textbf{Output}: ${\bm q}_k$\\
  \textbf{begin}: $l:=0$\\
  $\text{Screen-step}$: If $\pi_i < \epsilon$ ($i=1,\cdots,k$), we discard the corresponding model. Then \\
   $~~~~~~~~~~$ $k\rightarrow k_{\text{new}}$. For convenience, $k_{\text{new}}$ is also written as $k$. When $k=k_{min}$, we will
    $~~~~~~~~~~$ terminate  the screen-step.\\
  $\text{E-step}$: Compute $\bb J^l({\bm q}_k^{h_l},k)$ by (\ref{BYY}) and $(\gamma_{i,j})$ by (\ref{Estep}).\\
  $\text{M-step}$: Update the components of ${\bm q}_k^{h_l}$ by (\ref{Mstep}).\\
  $\text{S-step}$: Update the smoothing parameter $h_l$ by (\ref{smooth}). Calculate $\bb J^{l+1}({\bm q}_k^{h_{l+1}},k)$.\\
  $~~~~~~~~~~$ $\text{if}$ $\bb J^{l+1}({\bm q}_k^{h_{l+1}},k)<\bb J^l({\bm q}_k^{h_l},k)$ \\
  $~~~~~~~~~~~~$ $l=l+1$; else break;\\
   $~~~~~~~~~~$ $\text{end}$ \label{SmEM}
\end{algorithm}

\subsection{Post-processing for discrete structure}
For  discrete structures, the immediate  results by above  method  may be not good enough. To this end, we use  a post-processing based on  regularization  to improve the connectivity of the important features. It is implemented in a block-by-block manner, so we perform the post-processing for each gridblock. Let
${\widetilde A}_i$ denote the value of the $i$-th gridblock. The goal of post-processing is to minimize  the following function with respect to $A$,
\begin{equation}\label{post-pro}
\bb G(A)=({\widetilde A}_i-A)^{2}+\tau \bb T(A), \quad A\in\big[A^l,A^u\big],\quad 0<\tau<1,\quad i=1,\cdots,N_m,
\end{equation}
where $\bb T(A)$ denotes the regularization term, $\tau$ the corresponding regularization weight and $[A^l,A^u]$ is the domain of definition. We need to find the minimizer of
$\bb G(A)$, which is the estimation  for the $i$-th gridblock. For equation (\ref{post-pro}), the penalty term $\bb T(A)$ is applied to penalize values away from $A^l$ or $A^u$. Thus, the selection of the regularization term is important. In practice, we often use the following quadratic form for $\bb T(A)$,
\[
\bb T(A)=(b_1A-b_2)(b_3-b_4A),
\]
where $b_i$ $(i=1,2,3,4)$ depends on $A^l$ and $A^u$. Then, $\bb G(A)$ can be expressed as a quadratic form, i.e.,
\[
\bb G(A)=(1-\tau b_1b_4)\bigg(A-\frac{2{\widetilde A}_i-\tau(b_1b_3+b_2b_4)}{2(1-\tau b_1b_4)}\bigg)^2+c,
\]
where $c$ is a constant independent of $A$. To obtain  the global minimum, we impose  the constraint
$1-\tau b_1b_4>0.$
For the minimizer of the convex function $\bb G(A)$ in each gridblock, we note that there exist three cases:
\[
\begin{aligned}
\text{(i) If }&\frac{2{\widetilde A}_i-\tau(b_1b_3+b_2b_4)}{2(1-\tau b_1b_4)}<A^l,\quad {\widehat A}_i=A^l;\\
\text{(ii) If }&A^l\leq\frac{2{\widetilde A}_i-\tau(b_1b_3+b_2b_4)}{2(1-\tau b_1b_4)}\leq A^u,\quad{\widehat A}_i=\frac{2{\widetilde A}_i-\tau(b_1b_3+b_2b_4)}{2(1-\tau b_1b_4)};\\
\text{(iii) If }&\frac{2{\widetilde A}_i-\tau(b_1b_3+b_2b_4)}{2(1-\tau b_1b_4)}>A^u,\quad {\widehat A}_i=A^u,
\end{aligned}
\]
where $i=1,\cdots,N_m$.
Then, $\bb G({\widehat A}_i)$ is the global minimum in the $i$-th gridblock and ${\widehat A}_i$ is the estimation  in  the $i$-th gridblock.

\section{IES-IS  for  non-Gaussian priors}\label{IES-IS-non}
In this section, we present IES-IS for Bayesian inversion with priors described by  DCT and GMM.
 IES is a sampling method, which can generate ensemble samples to efficiently  estimate  the MAP point.
 We note that IES is  a Gaussian approximation.  To improve the effectiveness of ensemble samples, we use  IS method to get  a data-informed importance function. IS does not depend on any Gaussian assumption and is an importance sampling method, which can find the samples with high probability. The importance samples can be generated by  the  implicit equation, where the MAP point of $p_i(\bt|\bm d)$ and Hessian matrix of $\bb F_i(\bt)$ $(i=1,\cdots,k)$ are necessary. To avoid the computation of Jcoby matrix, we use the ensemble mean as the MAP point and approximate the inverse of Hessian matrix by Monte Carlo method. We perform a resampling method to avoid the ensemble degeneracy. The proposed IES-IS method is used to deal with the non-Gaussian Bayesian inverse problems through using  DCT and SmEM-based GMM.

\subsection{IES-IS based on DCT}
In this section, we use DCT  to parameterize the unknown function  $a(x)$ and obtain the prior information.  Each column of the basis function matrix $\bm\Phi_{N_m}$ is given by (\ref{DCT}). To perform the proposed algorithm efficiently, we use the truncated DCT expansion in (\ref{TDCT}).  The target field $a(x)$ is parameterized by $\bt$ through DCT.  Besides, the post-processing is applied to improve the connectivity of the inversion field. To further improve the efficiency, we use a criterion to reduce the dimension of $\bt$ against the iterations.    

Let  the  prior $p(\bt)$ be  Gaussian. Then this   corresponds to the case of $k=1$ in GMM described in Subsection \ref{IES-ISpr}.
 At the $l$-th iteration step, we use the Monte Carlo method in \cite{li2009iterative} to approximate the covariance matrix. Then 
 \begin{equation}\label{IES_DCT}
 \left\{
 \begin{aligned}
 &{\bm C}_{\theta_l}\approx \frac{1}{N_e}\sum_{j=1}^{N_e}({\bt}_l^j-\bar{{\bt}}_l)({\bt}_l^j-\bar{{\bt}}_l)^T\\
 &{\bm C}_{\theta_lD_l}\approx {\bm C}_{\theta_l}{\overline{\bm G}}_{l}^T\approx\frac{1}{N_e}\sum_{j=1}^{N_e}({\bt}_l^j-\bar{{\bt}}_l)\bigg(g(\widehat{\bm a}_l^j)-g(\bar{\bm a}_l)\bigg)^T\\
 &{\bm C}_{D_lD_l}\approx \overline{{\bm G}}_{l}{\bm C}_{\theta_l}\overline{{\bm G}}_{l}^T\approx\frac{1}{N_e}\sum_{j=1}^{N_e}\bigg(g(\widehat{\bm a}_l^j)-g(\bar{{\bm a}}_l)\bigg)\bigg(g(\widehat{\bm a}_l^j)-g(\bar{{\bm a}}_l)\bigg)^T,\\
\end{aligned}
\right.
 \end{equation}
where $\overline\bt_{l}=\frac{1}{N_e}\sum_{j=1}^{N_e}\bt_{l}^j$ and $\overline{\bm a}_l=\sum_{j=1}^{N_e}\widehat{\bm a}_l^j$.
By substituting (\ref{IES_DCT}) into (\ref{IES_kal}), (\ref{IES}) and (\ref{IS_H}),  we obtain the intermediate ensemble samples
\begin{equation}\label{IES_IS1}
\widetilde{\bt}_{l+1}^j={\bt}_l^j-\frac{1}{1+\lambda}({\bm C}_{\theta_l}-{\bm K}_l{\bm C}_{\theta_lD_l}^T){\bm C}_{\theta}^{-1}({\bt}_l^j-{\bt}^{pr})-{\bm K}_l\bigg(g(\widehat{\bm a}_l^j)-{\bm d}\bigg),\quad j=1,\cdots,N_e.
\end{equation}
The inverse of modified Hessian matrix can be approximated by
\[
\widetilde{{\bm H}}_{1}^{-1}=\frac{1}{1+\lambda}({\bm C}_{\theta_l}-{\bm K}_l{\bm C}_{\theta_lD_l}^T).
\]
The MAP point $\widetilde{\bm\mu}$ is obtained by equation (\ref{IS_mu}).
Thus we can obtain the importance ensemble samples $\{\bt_{l+1}^{j}\}_{j=1}^{N_e}$, which are generated by equation (\ref{sam_IS}) with weights obtained by equation(\ref{weight}).

To avoid the ensemble degeneracy, the resampling method is performed based on the  weights. Then we sort the components of ensemble mean $\overline{\bt}_{l+1}$ in descending order, i.e.,
\[
|\overline{\theta}^{(1)}_{l+1}|\geq|\overline{\theta}^{(2)}_{l+1}|\geq\cdots\geq|\overline{\theta}^{(N_l)}_{l+1}|,
\]
where $|\cdot|$ is the absolute operator.
Let $\alpha$ be a threshold value (e.g., $\alpha=0.95$).   We truncate $\overline{\bt}_{l+1}$ by the criterion
\begin{equation}\label{reduce}
\frac{\sum_{i=1}^{N_{l+1}}|\overline{\theta}^{(i)}_{l+1}|}{\sum_{i=1}^{N_l}|\overline{\theta}^{(i)}_{l+1}|}\geq \alpha, 
\end{equation}
where $N_{l+l}$ is the dimension of $\bt$ for the next iteration and $N_{l+1}\leq N_l$. Then we retain $N_{l+1}$ rows of ensemble samples according to truncated $\overline\bt_{l+1}$. For DCT basis functions, we select the corresponding columns to construct $\bm\Phi$ for the next iteration. For convenience, each truncated sample is still denoted by $\bt_{l+1}^j$.
Thus the dimension of the parameter decreases gradually in IES-IS algorithm. This can improve the computation efficiency.

Let $\bm A$ represent the field  $\ln a(x)$ in discrete sense.  To update the  field, we use 
 \begin{equation}\label{pervalue}
 \widetilde{\bm A}_{l+1}^j={\bm\Phi}_{l+1}{\bt}_{l+1}^j,\quad j=1,\cdots,N_e.
 \end{equation}
Due to the continuous values of $\widetilde{\bm A}_{l+1}^j$, we perform the post-processing to get discrete values and improve the connectivity of main features. For each coefficient sample,
\begin{equation}\label{IS_pos}
\widehat{\bm A}_{l+1}^j=\textbf{Post-processing}(\widetilde{\bm A}_{l+1}^j),\quad j=1,\cdots,N_e.
\end{equation}
The outline of IES-IS based on DCT is presented in Algorithm \ref{IES-IS-DCT}.
\begin{algorithm}
\caption{ DCT-based IES-IS algorithm}
  \textbf{Input}: The initial basis functions matrix ${\bm\Phi}_{N_0}$, prior ensemble $\{{\bt}_0^j\}_{j=1}^{N_e}$, $N_e$ the \\
  $~~~~~~~~~~~$ ensemble size, the initial $\lambda$, decay factor $\nu$, the prior mean ${\bt}^{pr}$, the prior\\
  $~~~~~~~~~~~$ covariance ${\bm C}_{\theta}$, sufficiently small $\epsilon$ and $N_0=N_c$.\\
  \textbf{Output}: final ensemble $\{{\bt}_l^j\}_{j=1}^{N_e}$ and
  $\{\widehat{\bm A}_{N_{l}}^j\}_{j=0}^{N_e}$\\
  $\bf1$. Compute the ${\widehat{\bm A}}_{N_0}^j={\bm\Phi}_{N_0}{\bt}_{0}^j,\quad j=1,\cdots,N_e$.\\
  \textbf{begin the iterations :} $l=0$\\
  $\bf2$. Calculate the intermediate ensemble samples $\{\widetilde{\bt}_{l+1}^j\}_{j=1}^{N_e}$ using equation (\ref{IES_IS1}).\\
  $\bf3$. Compute $\widetilde{\bm\mu}_{1}=\frac{1}{N_e}\sum_{j=1}^{N_e}\widetilde{\bt}_{l+1}^j$. Then draw samples $\bm\xi_j$ from the reference function\\
    $~~~~$ $B({\bm\xi})$ and derive the importance ensemble samples $\{\bt_{l+1}^j\}_{j=1}^{N_e}$ by equation (\ref{sam_IS}).\\
  $\bf4$. Update the weights $\{w_j\}_{j=1}^{N_e}$ using equation (\ref{weight}) and normalize the weights.\\
  $\bf5$. Resampling based on the weights and obtain a new ensemble. For convenience, \\
    $~~~~~$ the new importance ensemble is also denoted as $\{{\bt}_{l+1}^{j}\}_{j=1}^{N_e}$.\\
  $\bf6$. Calculate $\overline{{\bt}}_{l+1}$ and truncate $\overline{{\bt}}_{l+1}$ by equation (\ref{reduce}). Then $N_l\rightarrow N_{l+1}$.\\
  $\bf7$. Calculate $\{\widehat{\bm A}_{N_{l+1}}^j\}_{j=1}^{N_e}$ by equation (\ref{pervalue}) and (\ref{IS_pos}).\\
  $\bf8$. $\lambda=\frac{\lambda}{\nu}$ and $l:=l+1$.\\
  $\bf9$. Repeat $\bf2-\bf8$ until $\|\bar\bt_{l+1}-\bar\bt_l\|<\epsilon$.
  \label{IES-IS-DCT}
\end{algorithm}

%%%%%%%%%%%%%%%%%%%%%%%%%%%%%%%%%%%%%%%%%%%%%%%%%%%
\subsection{IES-IS based on GMM}
 Algorithm \ref{IES-IS-DCT} is based on the prior parameterized by  DCT.  To further extend  the proposed method to other  non-Gaussian models, we couple IES-IS with GMM.
We apply  SmEM algorithm to ensemble samples for estimating GMM parameters. Compared with conventional  IES method, the probability density function is the sum of $k$  Gaussian density functions at each iteration. The prior ensemble $\{{\bt}_0^j\}_{j=1}^{N_e}$ is drawn from Gaussian distribution, the support of which is $\bb R^{N_\theta}$. The initial parameter ${\bm q}=\{\pi_i,{\bm\mu}_i,{\bm\Sigma}_i\}_{i=1}^k$ is arbitrarily given. For each iteration in the proposed algorithm, we give a clustering step of $\bt$.
In this section, the initial $\bt^{pr,i}=\bt^{pr}$ and $\bm C_{\theta,i}=\bm C_{\theta}$.

At the $l$-th iteration step, SmEM method is applied to $\{{\bt}_l^j\}_{j=1}^{N_e}$ and  gives us the forecast
$\bm q^f=\{\pi_i^f,{\bm\mu}_i^f,{\bm\Sigma}_i^f\}_{i=1}^k$ and the membership probability matrix $(\gamma_{i,j})$, which can be expressed by
\[
\gamma_{i,j}=\frac{\pi_ip({\bt}_l^j|\bm\mu_i,\bm\Sigma_i)}{\sum_{m=1}^k\pi_mp({\bt}_l^j|\bm\mu_m,\bm\Sigma_m)},\quad i=1,\cdots,k,\quad j=1,\cdots,N_e.
\]
Then we can get a forecast $p^f({\bt})$ for  $p(\bt)$, which is a GMM, i.e.,
\begin{equation}\label{gmm}
p^f({\bt})=\sum_{i=1}^k\pi_i^fp({\bt}|{\bm\mu}_i^f,{\bm\Sigma}_i^f).
\end{equation}
The samples $\{{\bt}_l^j\}_{j=1}^{N_e}$ are  drawn from $p^f(\bt)$.
Let ${\overline{\bm G}}_{l,i}$ denote the Jacobian matrix at $\bm\mu_i^f$ and $n_i=\sum_{j=1}^{N_e}\gamma_{i,j}$. We use the Monte Carlo method to approximate the covariance matrix of each model $i$, i.e.,
\begin{equation}\label{CGMM}
\left\{
 \begin{aligned}
 {\bm C}_{\theta_l}^i &\approx\frac{\sum_{j=1}^{N_e}\gamma_{i,j}({\bt}_l^j-\bm\mu_i^f)({\bt}_l^j-\bm\mu_i^f)^T}{n_i}\\
{\bm C}_{\theta_lD_l}^i &\approx \bm C_{\theta_{l}}{\overline{\bm G}}_{l,i}^T\approx\frac{\sum_{j=1}^{N_e}\gamma_{i,j}({\bt}_{l}^j-\bm\mu_i^f)\bigg(g({\bt}_{l}^j)-g(\bm\mu_i^f)\bigg)^T}{n_i}\\
\bm C_{D_lD_l}^i &\approx {\overline{\bm G}}_{l,i}\bm C_{\theta_l}{\overline{\bm G}}_{l,i}^T\approx\frac{\sum_{j=1}^{N_e}\gamma_{i,j}\bigg(g({\bt}_{l}^j)-g(\bm\mu_i^f)\bigg)\bigg(g({\bt}_{l}^j)-g(\bm\mu_i^f)\bigg)^T}{n_i}
\end{aligned}
\right.
\end{equation}
and $\bm C_{D_l\theta_l}^i= (\bm C_{\theta_lD_l}^i)^T$.
By substituting  (\ref{CGMM}) into (\ref{IES_kal}), (\ref{IES}) and (\ref{IS_H}),
then intermediate ensemble samples for model $i$ can be obtained by
\begin{equation}\label{IESGMM}
\widetilde{\bt}_{l+1}^{j}={\bt}_{l}^j-\frac{1}{1+\lambda}({\bm C}_{\theta_l}^i-{\bm K}_l^i{\bm C}_{D_l\theta_l}^i){\bm C}_{\theta}^{-1}({\bt}_{l}^j-{\bt}^{pr})-{\bm K}_l^i\bigg(g({\bt}_{l}^j)-{\bm d}\bigg),\quad j=1,\cdots,N_e.
\end{equation}
Then the inverse of modified Hessian matrix can be approximated by
\[
\widetilde{{\bm H}_i}^{-1}=\frac{1}{1+\lambda}({\bm C}_{\theta_l}^i-{\bm K}_l{\bm C}_{D_l\theta_l}^i).
\]
For model $i$, we perform IS method and  resampling. The MAP point $\widetilde{\bm\mu}_i$ is obtained  by equation (\ref{IS_mu}). Thus the importance samples $\{\bt_{l+1}^{j,i}\}_{i=1}^{N_e}$ can be generated by equation (\ref{sam_IS}) with the weights $\{w_j^i\}$ obtained by equation (\ref{weight}).
The resampling based on the improved  weights can avoid ensemble degeneracy.  Using the membership probability matrix $(\gamma_{i,j})$, we combine $k$ ensembles together to form the posterior ensemble, i.e.,
\begin{equation}\label{combineGMM}
{\bt}_{l+1}^j=\sum_{i=1}^k\gamma_{i,j}{\bt}_{l+1}^{j,i}.
\end{equation}

In the analysis step of GMM, the means and covariances of the mixture models can be updated by
\begin{equation}\label{mcGMM}
\bm\mu_i=\frac{\sum_{j=1}^{N_e}\gamma_{i,j}{\bt}_{l+1}^{j,i}}{n_i},\quad\bm\Sigma_i=\frac{\sum_{j=1}^{N_e}\gamma_{i,j}({\bt}_{l+1}^{j,i}-\bm\mu_i)({\bt}_{l+1}^{j,i}-\bm\mu_i)^T}{n_i},\quad i=1,\cdots,k.
\end{equation}
Let $\overline{\bm G}_i$ denote the Jacobian matrix at $\bm\mu_i$. In the paper, we use the difference method to calculate $\overline{\bm G}_i$.
The weight of each model in (\ref{gmm}) can be updated based on the observation data ${\bm d}$, i.e.,
\begin{equation}\label{weGMM}
\pi_i=\frac{p({\bm d}|\bm\mu_i,\bm\Sigma_i)n_i}{\sum_{i=1}^kp({\bm d}|\bm\mu_i,\bm\Sigma_i)n_i},
\end{equation}
where
\[
p({\bm d}|\bm\mu_i,\bm\Sigma_i)=\frac{\exp{\bigg[-\frac{1}{2}\bigg({\bm d}-g(\bm\mu_i)\bigg)^T({\overline{\bm G}}_i\bm\Sigma_i\overline{\bm G}_i^T+{\bm C}_D)^{-1}\bigg({\bm d}-g(\bm\mu_i)\bigg)\bigg]}}{\sqrt{(2\pi)^{N_d}\det(\overline{\bm G}_i\bm\Sigma_i\overline{\bm G}_i^T+\bm C_D)}}.
\]
Here, $\det(\cdot)$ denotes the determinant operator of a matrix.
When GMM-based IES-IS  achieves certain convergence, we obtain  a point estimate from the Gaussian mixture distribution,
 \begin{equation}\label{unknown_me}
 \widehat{\bt}=\mathbb{E}[{\bt}]=\sum_{i=1}^k\pi_i\bm\mu_i.
 \end{equation}
 The final posterior can be approximated by a GMM, i.e.,
 \[
 p(\bt|\bm d)\propto\sum_{i=1}^k\pi_ip(\bt|\bm\mu_i,\bm\Sigma_i).
 \]
 The pseudo-code for GMM-based IES-IS algorithm is provided in Algorithm \ref{IES-IS-GMM}.
\begin{algorithm}
\caption{ GMM-based IES-IS algorithm}
  \textbf{Input}: prior ensemble $\{{\bt}_0^j\}_{j=1}^{N_e}$, $N_e$ the ensemble size, the initial $\bm q=\{\pi_i,\bm\mu_i,\bm\Sigma_i\}_{i=1}^k$, \\
  $~~~~~~~~~$ the number of iterations $n$, sufficiently small $\epsilon$, positive integer $k_{max}$ and $k_{min}$, \\
  $~~~~~~~~~$ the initial $\lambda$, decay factor $\nu$, the prior mean ${\bt}^{pr}$ and the covariance $\bm C_{\theta}$.\\
  \textbf{Output}: final posterior ensemble $\{{\bt}_l^j\}_{j=1}^{N_e}$\\
  $l:=0$\\
  \textbf{begin the iterations :}\\
  $\bf1$. Perform the Algorithm \ref{SmEM} for $\{{\bt}_{l}^j\}_{j=1}^{N_e}$ to get $\bm q^f$ and $k$.\\
  $\bf2$. \textbf{for} $i=1:k$\\
    (1) Calculate the intermediate ensemble samples $\{\widetilde{\bt}_{l+1}^j\}_{j=1}^{N_e}$ using equation (\ref{IESGMM}).\\
    (2) Compute ${\widetilde{\bm\mu}}_{i}=\frac{1}{N_e}\sum_{j=1}^{N_e}\widetilde{\bt}_{l+1}^j$. Then draw samples $\bm\xi_j$ from the reference function\\
    $~~~~~$ $B({\bm\xi})$ and derive the importance ensemble samples $\{\bt_{l+1}^{j,i}\}_{j=1}^{N_e}$ by equation (\ref{sam_IS}).\\
    (3) Update the weights $\{w_j^i\}_{j=1}^{N_e}$ using equation (\ref{weight}) and normalize the weights.\\
    (4) Resampling based on the weights and obtain a new ensemble. For convenience, \\
    $~~~~~$ the new importance ensemble is also denoted as $\{{\bt}_{l+1}^{j,i}\}_{j=1}^{N_e}$.\\
    $~~~~$\textbf{end}\\
  $\bf3$. Update the posterior ensemble using equation (\ref{combineGMM}).\\
  $\bf4$. Update $\{\bm\mu_i,\bm\Sigma_i,\pi_i\}$ for the mixture posterior using equations (\ref{mcGMM})-(\ref{weGMM}).\\
  $\bf5$. $\lambda=\frac{\lambda}{\nu}$ and $l=l+1$.\\
  $\bf6$. Repeat $\bf1-\bf5$ until $\|\overline\bt_{l+1}-\overline\bt_{l}\|<\epsilon$.
  \label{IES-IS-GMM}
\end{algorithm}

Next, we present a result about GMM of the posterior when the prior is a GMM.
\begin{thm}
Assume that the forward operator $g(\cdot)$ is linear, and observation noise $\bm\varepsilon\sim{N}(\bm 0, \bm I)$. If  $\bm G=g$
and  the prior $p(\bt)$ is given by
\[
p(\bt)=\sum_{i=1}^k\pi_ip({\bt}|{\bm\mu}_i,{\bm\Sigma}_i),
\]
 then the posterior can be expressed by
\[
p(\bt|\bm d)\propto\sum_{i=1}^k\pi_i^ap({\bt}|{\bm\mu}_i^a,{\bm\Sigma}_i^a),
\]
where
\begin{equation}\label{thmGMM}
\left\{
 \begin{aligned}
{\bm\mu}_i^a &= (\bm G^T\bm C_D^{-1}\bm G+\bm\Sigma_i^{-1})^{-1}(\bm G^T\bm C_D^{-1}\bm d+\bm\Sigma_i^{-1}\bm\mu_i)\\
{\bm\Sigma}_i^a &= (\bm G^T\bm C_D^{-1}\bm G+\bm\Sigma_i^{-1})^{-1}\\
{\pi}_i^a &= \frac{\widetilde\pi_i}{\sum_{i=1}^{k}\widetilde\pi_i},
\end{aligned}
\right.
\end{equation}
and
\[
\widetilde\pi_i=\frac{\pi_idet(\bm\Sigma_i^a)^{-\frac{1}{2}}\exp\bigg[-\frac{1}{2}\bigg(\bm d^T\bm C_D^{-1}\bm d+\bm\mu_i^T\bm\Sigma_i^{-1}\bm\mu_i-(\bm\mu_i^a)^{T}(\bm\Sigma_i^a)^{-1}\bm\mu_i^a\bigg)\bigg]}{det(\bm\Sigma_i)^{\frac{1}{2}}det(\bm C_D)^{\frac{1}{2}}}.
\]
\begin{proof}
 Due to $\bm\varepsilon\sim{N}(\bm 0, \bm I)$,  we have the likelihood function
\[
p(\bm d|\bt)=\frac{1}{(2\pi)^{\frac{N_d}{2}}det(\bm C_D)^{\frac{1}{2}}}\exp\bigg[-\frac{1}{2}(\bm d-\bm G\bt)^T\bm C_D^{-1}(\bm d-\bm G\bt)\bigg].
\]
The prior $p(\bt)$ is a mixture of $k$ Gaussian densities. Using  Bayes rule, we get the posterior
\[
p(\bt|\bm d)\propto p(\bt)p(\bm d|\bt)=\bigg(\sum_{i=1}^k\pi_ip({\bt}|{\bm\mu}_i,{\bm\Sigma}_i)\bigg)p(\bm d|\bt)=\sum_{i=1}^k\pi_i\bigg(p({\bt}|{\bm\mu}_i,{\bm\Sigma}_i)p(\bm d|\bt)\bigg),
\]
where
\[
p(\bt|\bm\mu_i,\bm\Sigma_i)=\frac{1}{(2\pi)^{\frac{N_\theta}{2}}det(\bm\Sigma_i)
^{\frac{1}{2}}}\exp\bigg[-\frac{1}{2}(\bt-\bm\mu_i)^T\bm\Sigma_i^{-1}(\bt-\bm\mu_i)\bigg].
\]

We complete the perfect square of $\log \bigg(p({\bt}|{\bm\mu}_i,{\bm\Sigma}_i)p(\bm d|\bt)\bigg)$ with respect to $\bt$. In  the end, we  have the posterior with a GMM, the means, covariances and weights, which are   given by equation (\ref{thmGMM}). We note that  the weights have been normalized here.
\end{proof}
\end{thm}

When $g$ is nonlinear, we can not get a close  expression for  the posterior. In this case, we  use ensemble method to get a Gaussian approximation for each model $i$ and  obtain an approximation GMM of the posterior. The means, covariances and weights are given by equation (\ref{mcGMM}) and (\ref{weGMM}). The model number  $k$  is  given by BYY-HDS method.

%%%%%%%%%%%%%%%%%%%%%%%%%%%%%%%%%%%%%%%%%%%%%%%%%%%%%%%%%%%%%
%%%%%%%%%%%%%%%%%%%%%%%%%%%%%%%%%%%%%%%%%%%%%%%%%%%%%%%%%%%%%%%%

\section{Numerical examples}\label{exam}
In this section, we apply the proposed IES-IS to  subsurface flows  and anomalous diffusion problems in porous media and estimate the model's unknown inputs.
 In Subsection \ref{source}, we estimate the source locations of the single-phase flow using GMM-based IES-IS method. In Subsection \ref{channel}, we recover a channel structure in a permeability field by DCT-based IES-IS method. In Subsection \ref{fracture}, we will identify the fracture in porous media through GMM-based IES-IS method.

For the numerical examples,  we consider a dimensionless  square domain $\Omega=[0,1]\times[0,1]$ for spatial variable. Observation data are generated synthetically by using FEM in a fine time division and the forward problems are  solved by FEM with a coarse time division to avoid inverse crime.  The ensemble size $N_e$ is set as $2000$. The decay factor $\nu$ is set as 2 and the initial $\lambda$ is set as 1 in Algorithm \ref{IES-IS-DCT} and \ref{IES-IS-GMM}. The mean $\bt^{pr}=\bm0$ and the covariance $\bm C_\theta=\bm I$. In these examples, $I$ denotes  the number of iterations and the covariance matrix of observation error
$\bm C_D=\sigma^2\bm I$.

%%%%%%%%%%%%
\subsection{Estimate  source locations}\label{source}

In this subsection, we consider a steady single-phase flow model
\[
-\nabla\cdot(a(x)\nabla u(x))=f(x),\quad \quad x\in\Omega,
\]
with a mixed boundary condition, where Dirichlet boundary condition is
\begin{equation*}
u(0,y)=1 ,\quad \quad \quad u(1,y)=0,
\end{equation*}
and no flow boundary condition is imposed on the other two boundaries.  Here, the permeability field $a(x)$ is given by
$a(x)=\exp(1+0.5x+y).$
The source term has the form
\[
  f(x,\bt)=\frac{s}{\pi\iota}\exp\{-\frac{\|{\bt}-x\|^2}{2\iota^2}\},
\]
where $s=\exp(2)$ is the strength, $\iota=0.05$ is the width and ${\bt}$ is the unknown  source location.

\begin{figure}[tbp]
  \centering
  \includegraphics[width=5in, height=1.3in]{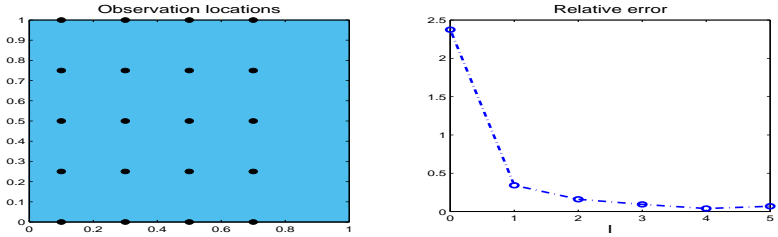}
  \caption{ The observation locations (left) and the relative error $\varepsilon_{\theta}$ (right). $I=0$ denotes the error obtained by the prior ensemble.}\label{perob}
\end{figure}
The truth source location is set as $\bt^{tr}=(0.09,0.23)^T$.
Observations are taken from the single-phase flow model, where the locations are distributed on the uniform $4\times 5$ grid of the domain $[0.1, 0.7]\times [0, 1]$ as shown in Figure \ref{perob} (left).
The forward model is solved on a uniform $100\times100$ grid and the observation data are obtained by solving the problem on  a uniform $200\times200$ grid.
The standard deviation $\sigma$ is set as 0.01 and the scale parameter $\rho=1$ in equation(\ref{weight}).
 GMM is used to characterize  the prior.   The initial $\{\bm\mu_i,\bm\Sigma_i\}_{i=1}^k$ is arbitrarily given and $\pi_i=\frac{1}{k}$. The initial $k$ is set as $k_{max}$. We take $k_{min}=2$ and $k_{max}=5$.

To measure the estimate accuracy, we  define the relative errors $\varepsilon_{\theta}$  by
\begin{equation*}
  \quad \varepsilon_{\theta}:= \frac{\|\widehat{\bt}-{\bt^{tr}}\|}{\|{\bt^{tr}}\|},
\end{equation*}
where $\widehat{\bt}$ is the estimation  given by equation (\ref{unknown_me}).  The mixture models are screened by SmEM method presented in Algorithm \ref{SmEM}, where the models with the relatively small weight will be abandoned. Thus, the model number  $k$ may decrease via the iteration. We  find the relative error of the prior ensemble is very large by Figure \ref{perob} (right).
When the observation data have been incorporated into the prior ensemble, the error decreases as more iterations are implemented.
  Then the relative error gradually tends to be stable against the iterations.  This demonstrates that Algorithm \ref{IES-IS-GMM} is effective and convergent.

\begin{figure}[tbp]
  \centering
  \includegraphics[width=5in, height=1.3in]{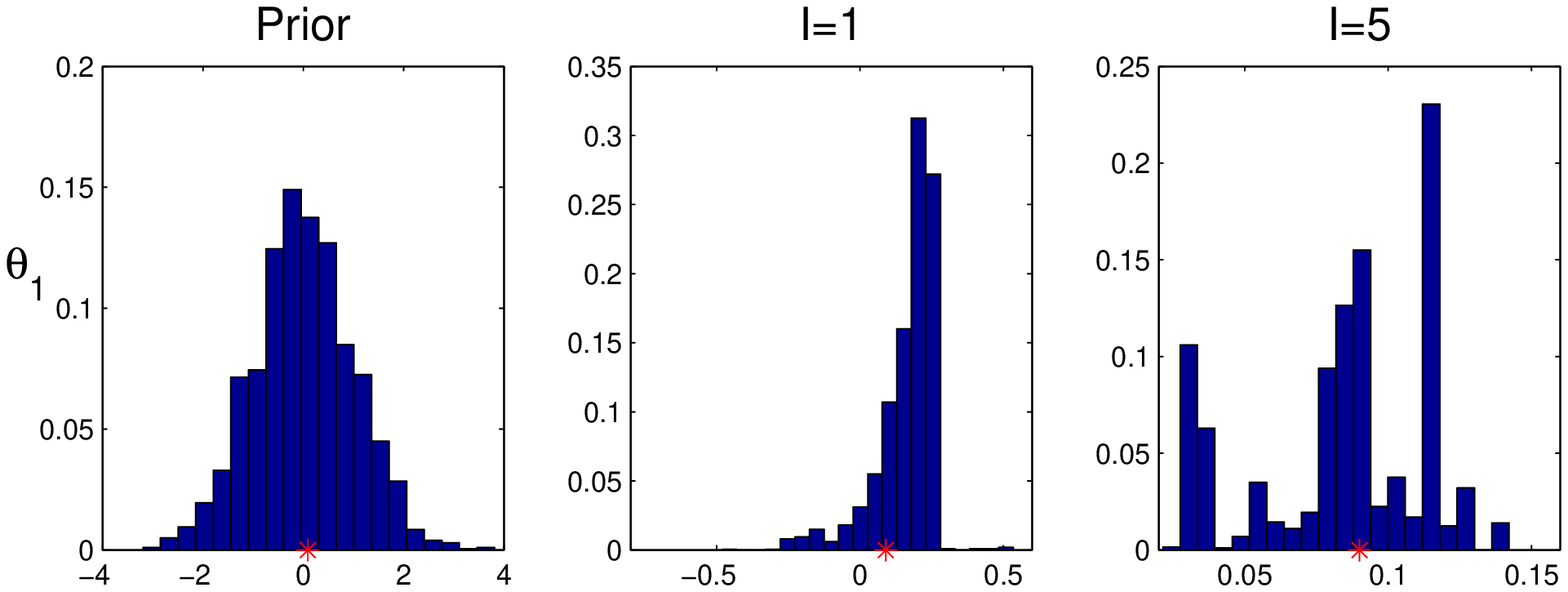}
  \includegraphics[width=5in, height=1.3in]{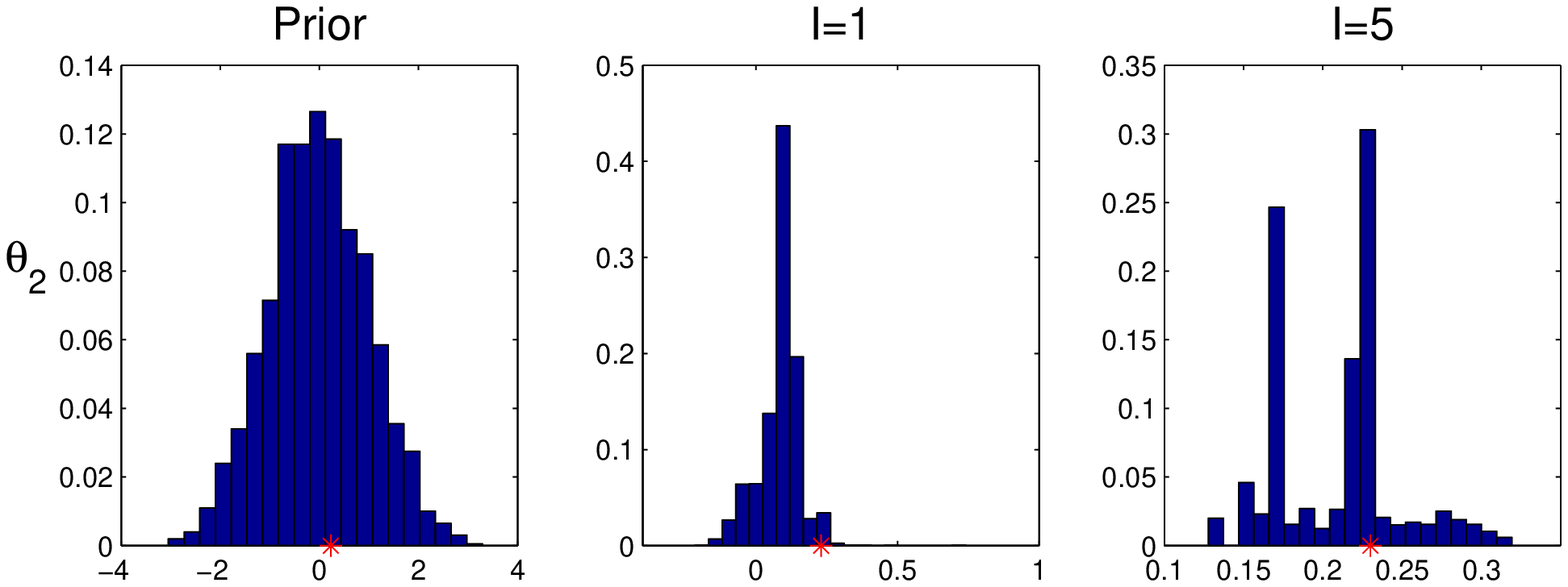}
  \caption{ The hists of marginal posterior density for $\bt$ at different iterations and the red markers are the reference values.}\label{density}
 \end{figure}
In GMM-based IES-IS algorithm, the posterior distribution is non-Gaussian. To account for the non-Gaussain property of the unknowns at different iterative steps, we plot the hists of marginal posterior density function in Figure \ref{density}.
The prior ensemble is randomly drawn from the Gaussian distribution. After the first iteration, the ensemble posterior becomes  non-Gaussian as shown in second column of Figure \ref{density}. For the first component of $\bt$, the skewness is obvious. This due to the importance ensemble samples are screened by the implicit sampling. We note that most samples cluster in one interval with large weights. Finally, the samples cluster in several intervals.
This shows the  non-Gaussian property  of the posterior ensemble.

\begin{figure}[tbp]
  \centering
  \includegraphics[width=5in, height=1.3in]{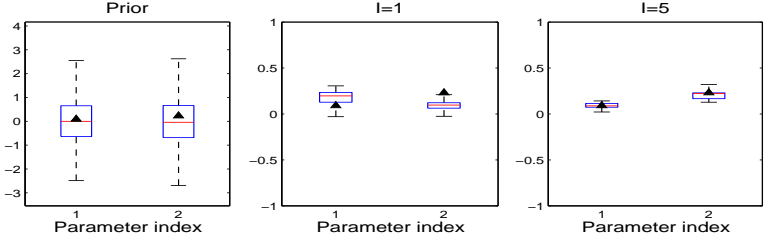}
  \caption{Medians (red solid line), reference (black filled triangle), the edges of the box (blue unfilled rectangle) are the 25th (bottom) and 75th (top) percentiles, and $95\%$ credible intervals (black dash line) for $\bt$ via the iteration.}\label{spreu}
\end{figure}
Figure \ref{spreu} shows the medians, percentiles and  $95\%$ credible intervals for parameter $\bt$ by the posterior ensemble samplers at different iterations.
We  find that  the medians of the prior ensemble is in the middle of the blue rectangle and the credible intervals are big at this stage. This is because  the samplers are drawn from the Gaussian distribution. After performing the proposed IES-IS, the medians deviate away from the middle of the rectangle. This implies that the posterior distributions of the unknowns are skewed and non-Gaussian. The credible intervals gradually become narrow when more iterations are used and the uncertainty of the unknowns decreases.
The reference values are included in the $95\%$ credible intervals. The marginal posterior densities are  skewed and  non-Gaussian.

\begin{figure}[tbp]
  \centering
  \includegraphics[width=5in, height=1.3in]{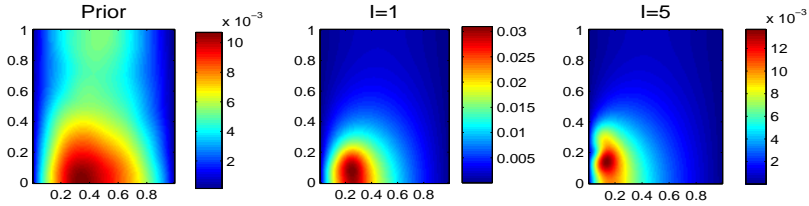}
  \caption{The standard deviation of the simulated states obtained by the posterior ensemble against the iteration steps.}\label{std}
\end{figure}
Figure \ref{std} depicts the standard deviation of the simulated state $u$ with the posterior ensemble samplers and realizations constructed by GMM-based IES-IS method.
For the prior ensemble samplers, the standard deviation is large and scattered in all physical domain except for the constrained boundaries. Then the area of large standard deviation becomes smaller  because the uncertainty of $\bt$ decreases. We note that the uncertainty of source location is close to $x=0$ when the iteration moves on. Finally, the uncertainty is mainly  concentrated in a small region  around the true source location. The standard deviation obtained by the final posterior ensemble is smallest.  This implies that the source location is identified.

\begin{figure}[tbp]
  \centering
  \includegraphics[width=5in, height=1.5in]{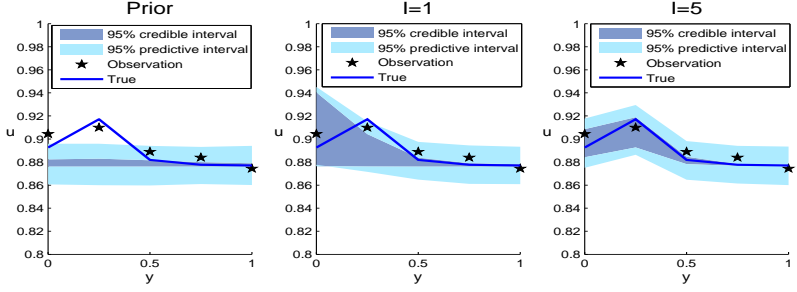}
  \caption{95\% prediction intervals, 95\% credible intervals, observations and true values for $u(0.1,y)$ at different iterations.}\label{ssim2}
\end{figure}
The credible and prediction intervals, along with the truth and observation data, are illustrated for $u(0.1,y)$ in Figure \ref{ssim2}.
We note that some observations are outside of the prediction intervals obtained by the prior ensemble,  which take account of the observational noise. When the information is gained  from the observation data, the shape of prediction and credible intervals for $u(0.1,y)$ align with  the true values and observation data. The credible interval of $u(0.1,y)$ becomes wider  as $x$ gets closer to $x=0$. This is due to the true source location is close to $x=0$. Finally, the true values and observation data are contained in the prediction intervals and credible intervals.

%%%%%%%%%%%%%%%%%%%%%%%%%%%%%%%%%%%%%%%%%%%%%%%%%%%%%%%%%%%%%%%%%%%%%%%%%%%%%%%%%%%%%%%%%%%%%%%%%%%%%%%%%%%%%%%%%%%%%%%%%%%%%%%%
\subsection{Recover  the channels structure}\label{channel}
\begin{figure}[tbp]
  \centering
    \includegraphics[width=1.6in, height=1.4in]{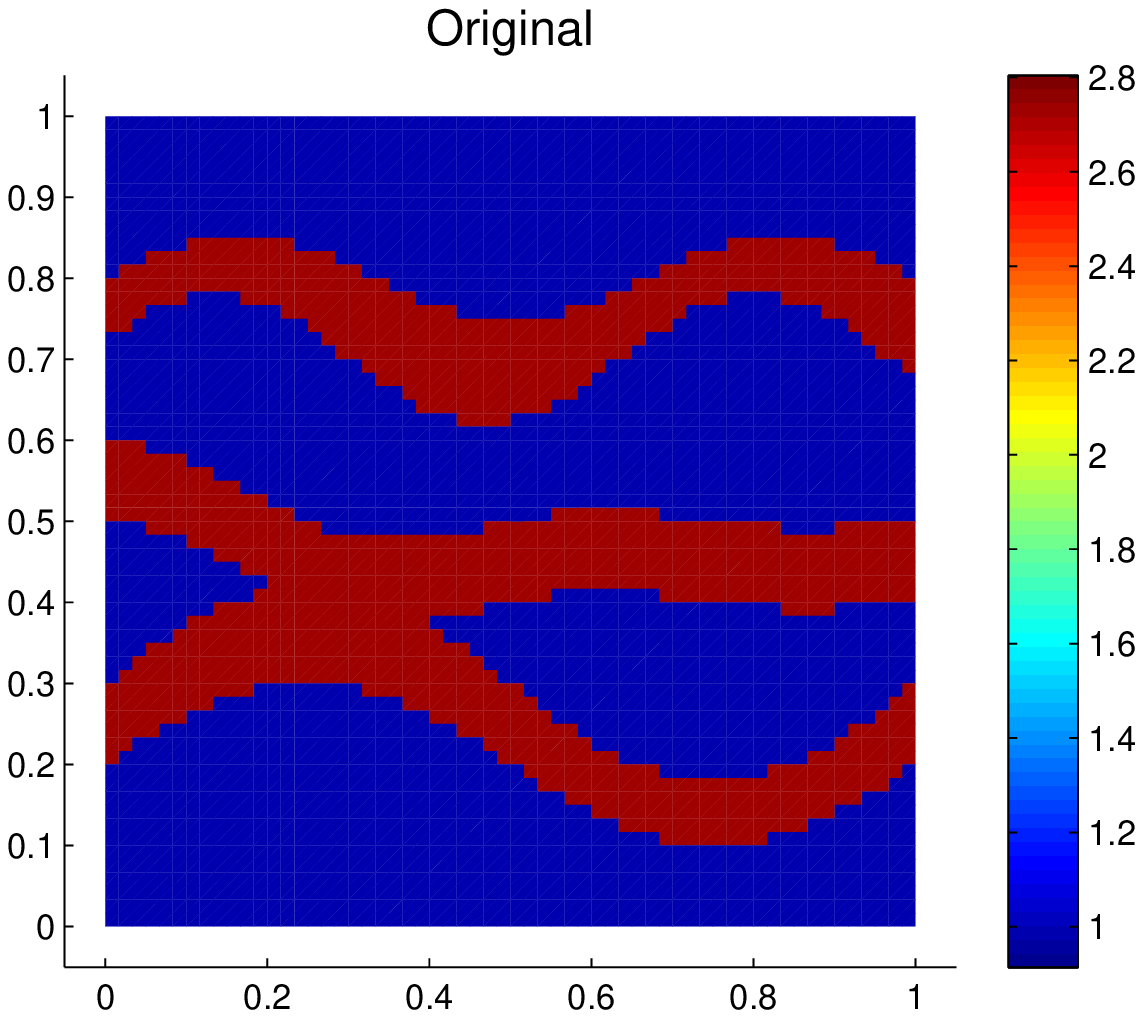}
  \includegraphics[width=1.6in,height=1.4in]{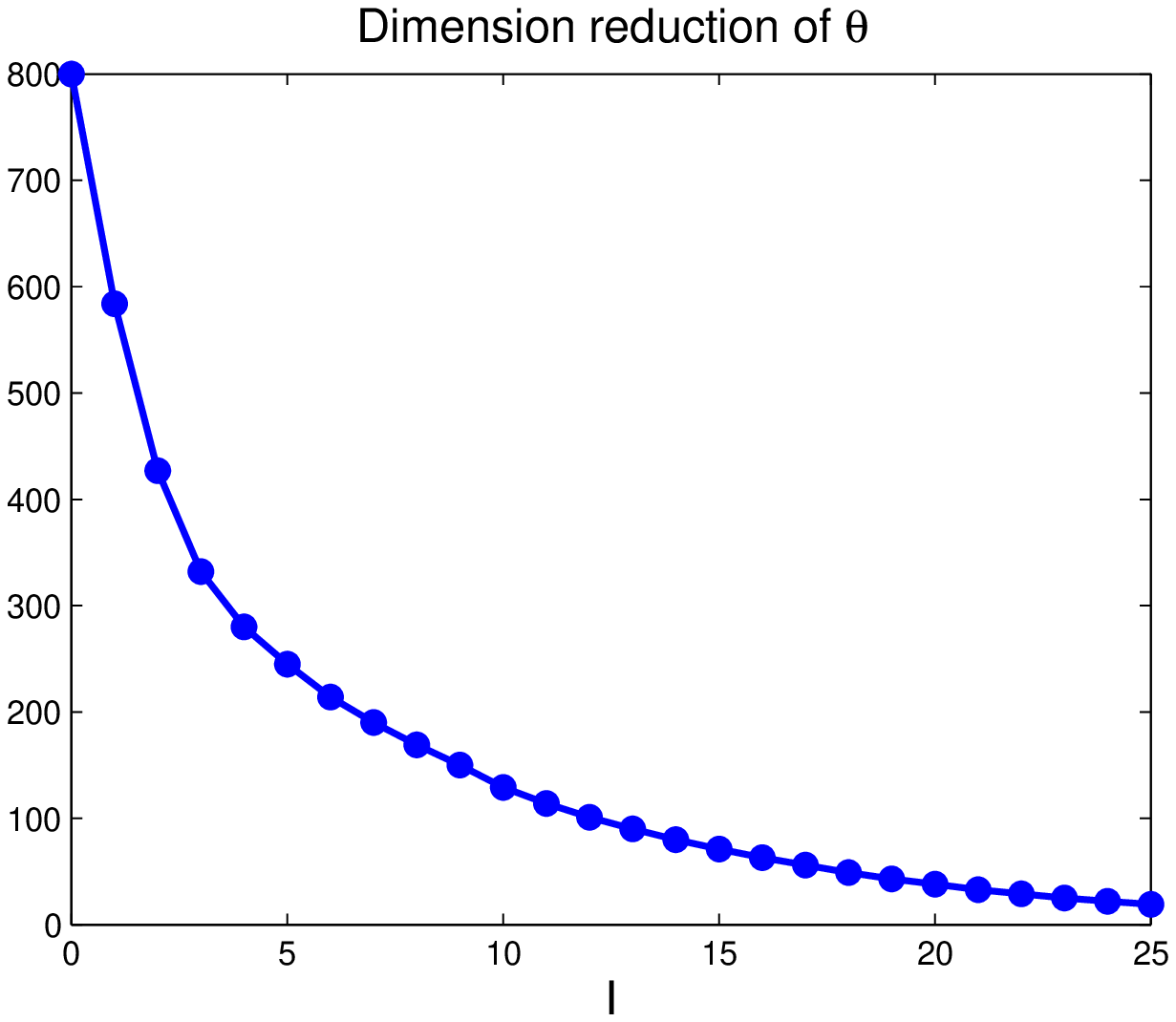}
  \includegraphics[width=1.6in,height=1.4in]{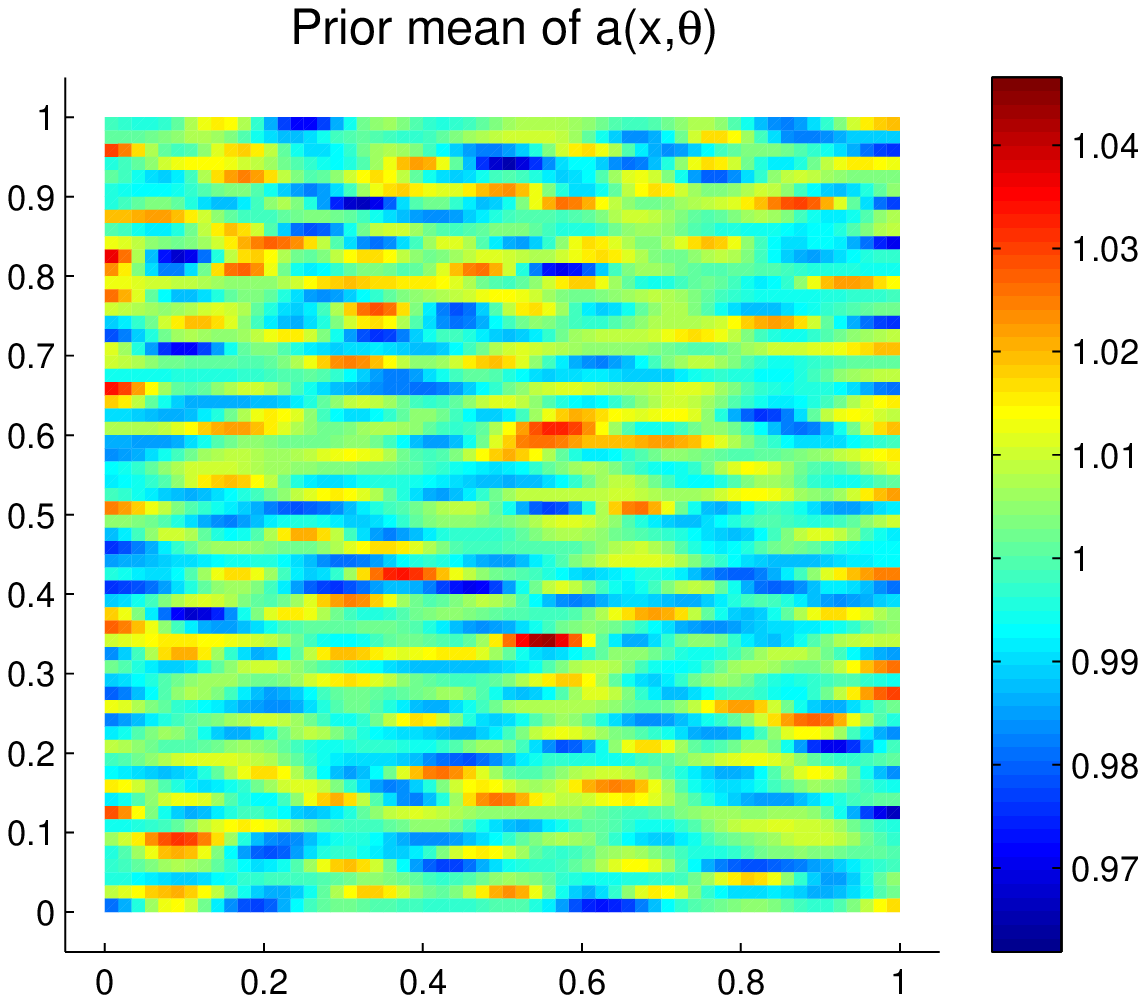}
  \caption{ The true profile of the permeability field $a(x)$ (left), the dimension reduction of $\bt$ (middle) and the prior mean of the permeability field $a(x,\bt)$ (right).}\label{chan}
\end{figure}

In this subsection, we consider the unsteady sing-phase  flow model
\[
\frac{\partial u(x,t)}{\partial t}-\nabla\cdot(a(x)\nabla u(x,t))=f(x,t),\quad \quad x\in\Omega, \quad t\in (0, T]
\]
with mixed boundary condition, where Dirichlet boundary condition is
\begin{equation*}
u(0,y;t)=1 ,\quad \quad \quad u(1,y;t)=0,
\end{equation*}
and there is no flow on the other boundaries. The source term is  $f=10$ and $T=1$. The permeability field $a(x)$ is unknown and needs to be recovered.
We   have the prior information of the permeability field, which is structured with the channels that lie between $y=0$ and $y=1$. The truth permeability field is divided into $7$ parts as shown in Figure \ref{chan} (left).
The forward model is solved by a uniform $60\times60$ grid with time step $\Delta t=0.02$, and the observations are obtained  with time step  $\Delta t=0.01$ to avoid the inverse crime. Here, observation data are taken at time instance $t=1$. We take $504$ observations to perform DCT-based IES-IS method. The standard deviation $\sigma$ is set as 0.01 and the scale $\rho=10$ in equation(\ref{weight}). The weight $\tau$ is set as 0.75 and the regularization terms in equation (\ref{post-pro}) is given by
\[
\bb T(A)=A(1-A),\quad A\in[0,1].
\]
To parameterize the channels, we use the truncated DCT expansion by the first 800 terms. Thus the logarithmic permeability field can be expressed as
\[
\log a(x)\approx \log a(x,\bt)=\sum_{r=0}^{N_c-1}\theta_r{\bm\phi}_r,\quad N_c=800.
\]
Thus we have the unknown parameter $\bt=(\theta_0,\cdots,\theta_{799})^T$.

To identify the channels structure, we combine IES-IS with DCT in this example. The dimension of $\bt$ for the original discrete model is 3600. Due to the large uncertainty of prior ensemble, we need enough basis functions to retain the main features of the channels.  Thus we select the first 800 columns of $\bm\Phi_{N_m}$ as the prior basis functions, which contain the important information of the channels. In IES-IS, we dynamically  reduce the dimension of $\bt$ by the posterior ensemble at each iteration, as shown in Figure \ref{chan} (middle). The dimension of parameters decreases against the iterations. When the uncertainty of the channels decreases, the reduction of the dimension becomes slow. Finally, we only use $19$ parameters  to construct the channels structure.

\begin{figure}
  \centering
  \includegraphics[width=1.6in,height=1.4in]{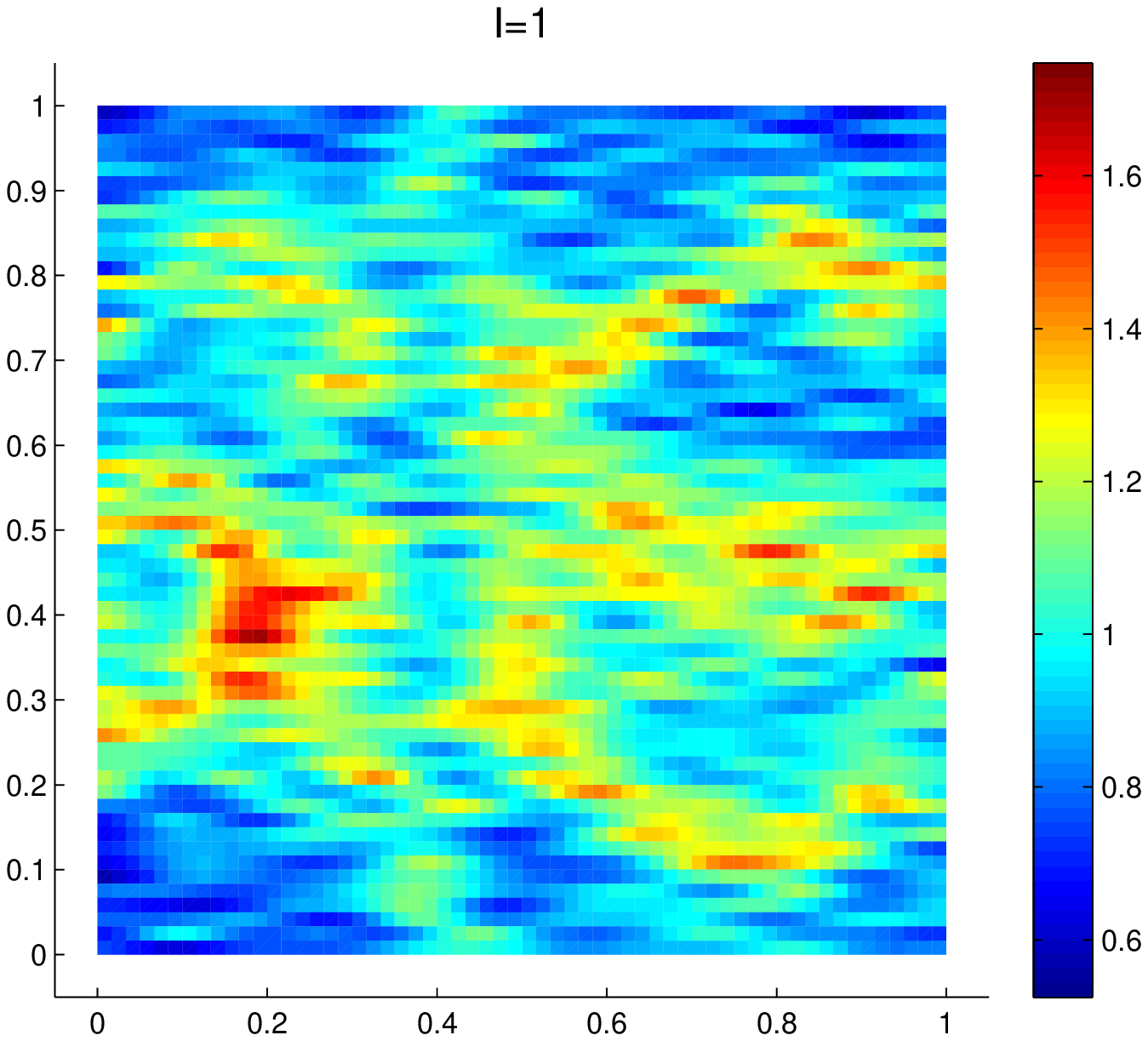}
  \includegraphics[width=1.6in,height=1.4in]{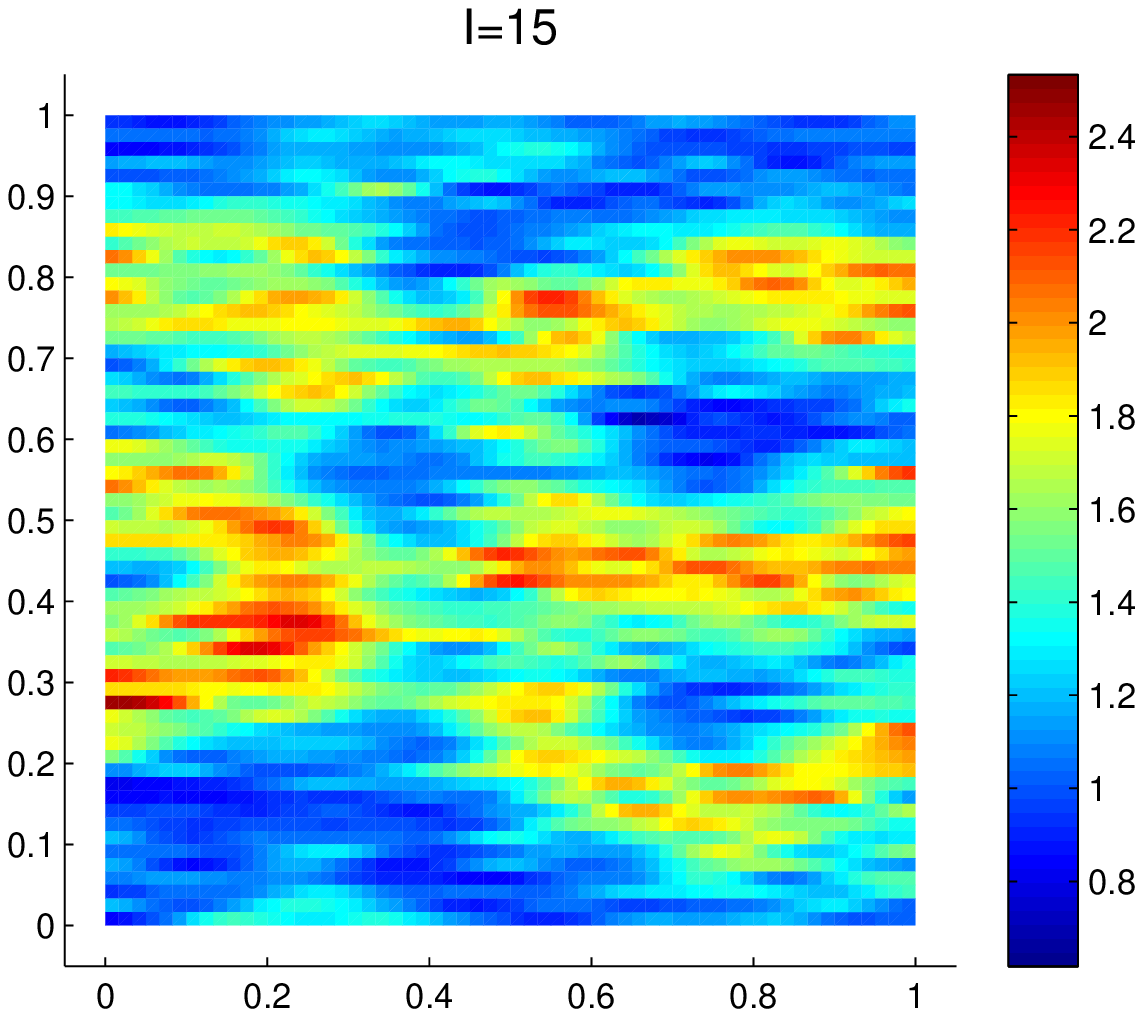}
  \includegraphics[width=1.6in,height=1.4in]{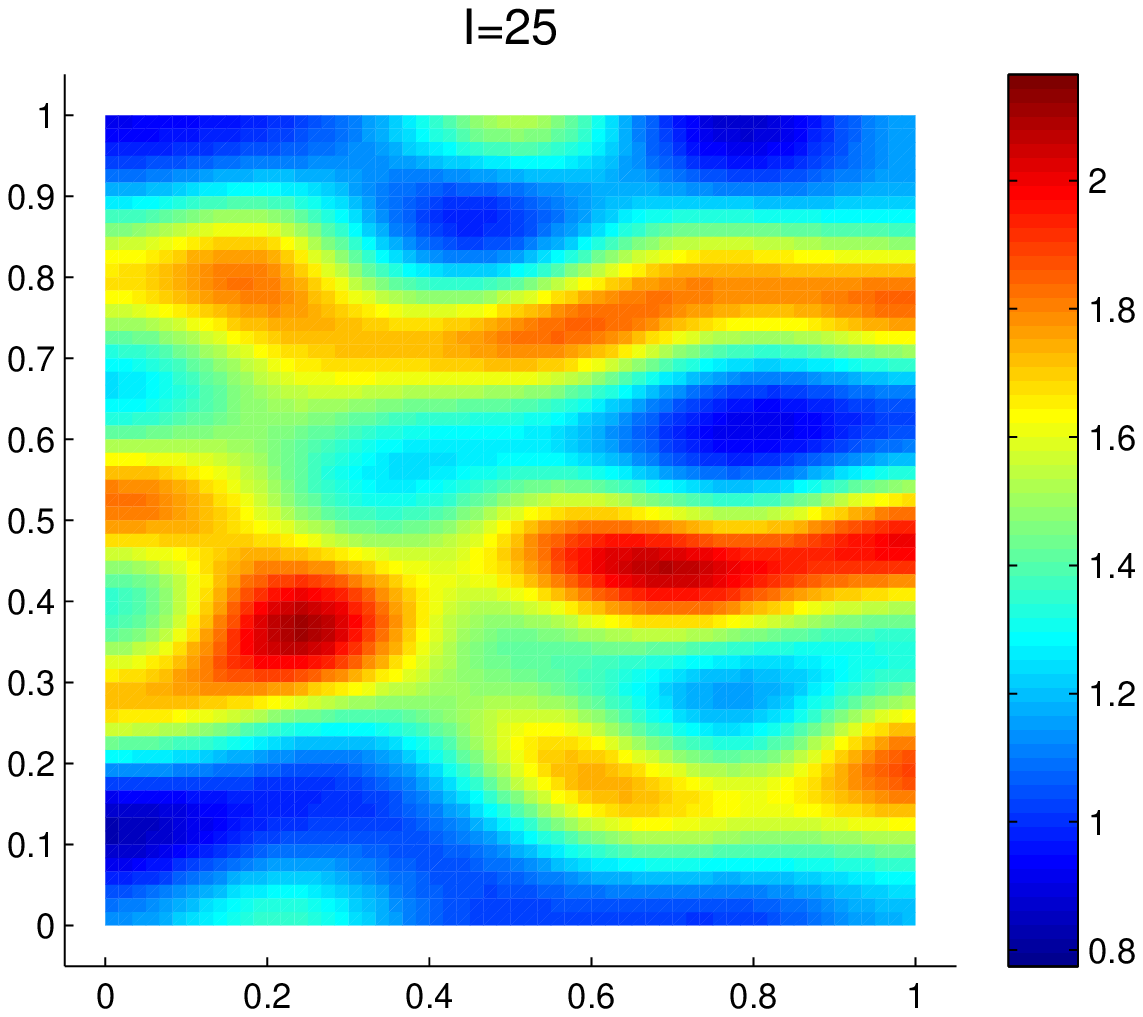}
  \includegraphics[width=1.6in,height=1.4in]{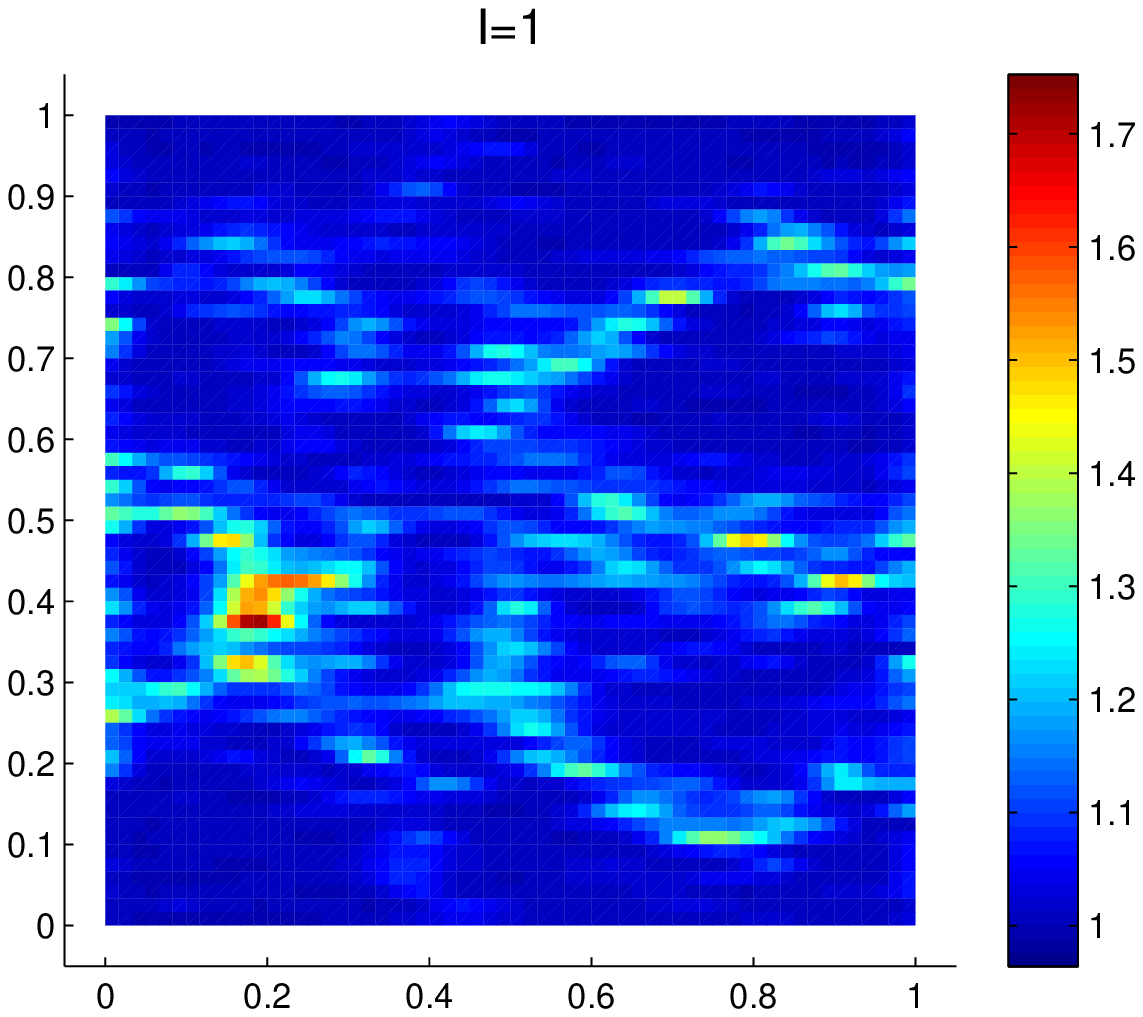}
  \includegraphics[width=1.6in,height=1.4in]{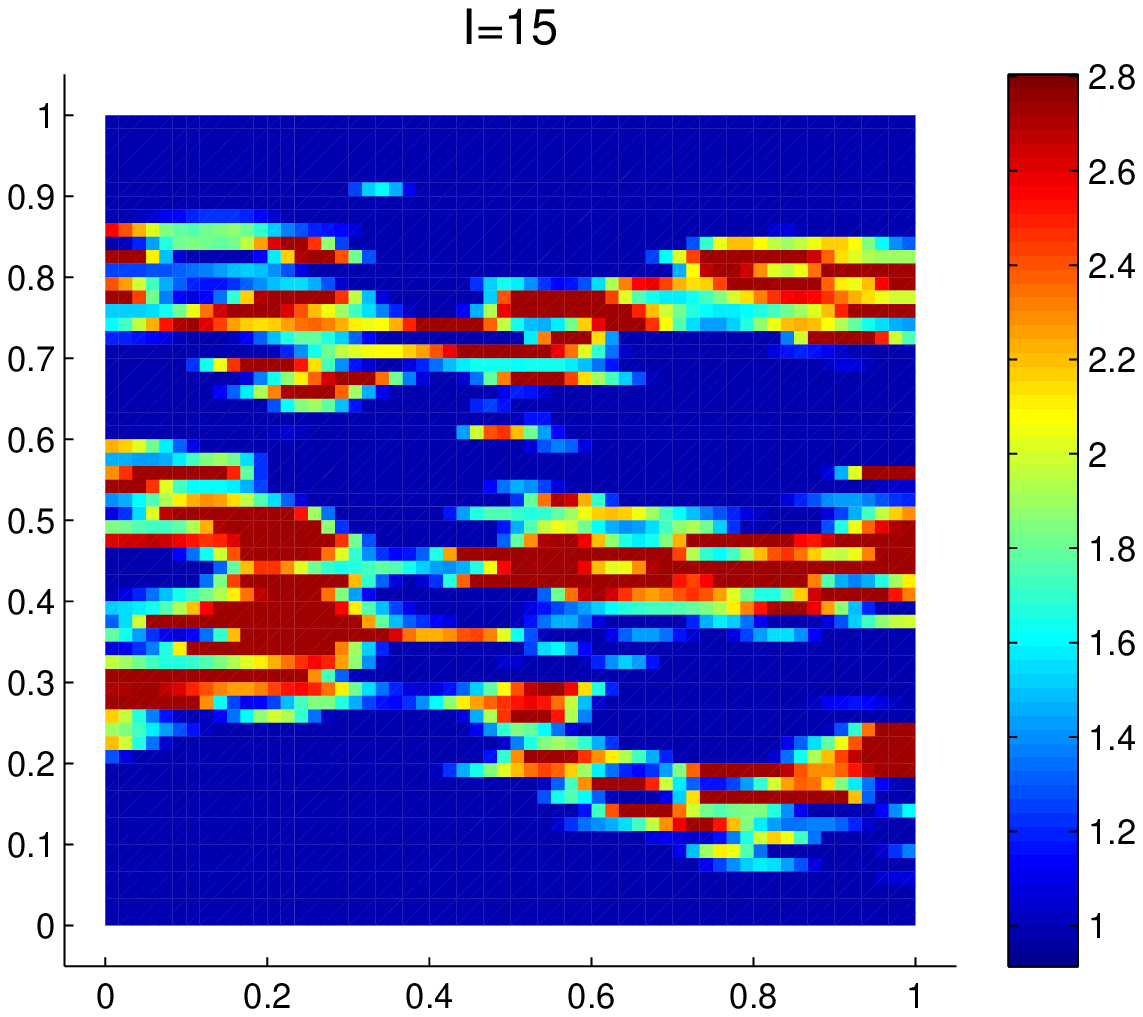}
  \includegraphics[width=1.6in,height=1.4in]{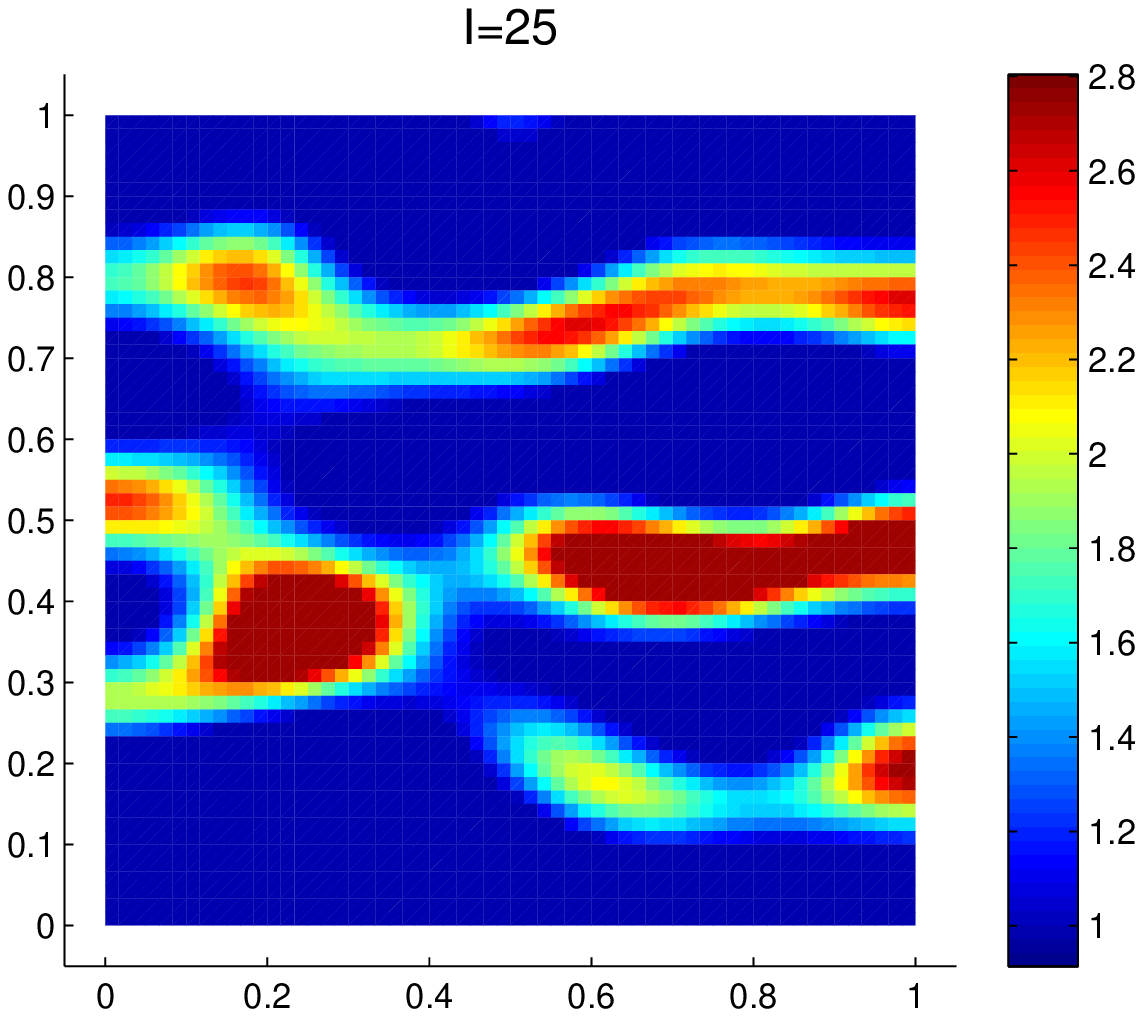}
  \caption{The first row is the posterior means of the permeability field $a(x,\bt)$ without the post-processing via the iterations. The posterior means of the permeability field $a(x,\bt)$ by the post-processing via the iterations in the second row.}\label{cha}
\end{figure}

The posterior means of the permeability field $a(x,\bt)$ via the iterations are presented in Figure \ref{cha},
The first row of which  illustrates the posterior mean without the post-processing, where the inversion  $a(x,\bt)$ looks  continuous. We see  that the prior ensemble does not give any channel structure shown in Figure \ref{chan} (right). When the information of observation data have been incorporated into the prior ensemble,  a ambiguous channel structure appears. Moreover, we  see that a  profile of the channels and the values in the channels get closer to the truth as more iterations are used. 
In order to show the effect  of the post-processing, we plot the posterior means of $a(x,\bt)$ using  the post-processing for the posterior ensemble in the second row of Figure \ref{cha}.
Compared with the first row, the values of input filed $a(x,\bt)$ with the post-processing are discontinuous. This is due to that the post-processing adopts the regularization method to reduce the continuity.  Compared with the first row,  the final channels are closer to the truth  channels in Figure \ref{chan} (left) and the values of permeability field are more accurate.  This shows that the post-processing based on a regularization  improve the construction of  the channels.

\begin{figure}
  \centering
  \includegraphics[width=5in,height=1.6in]{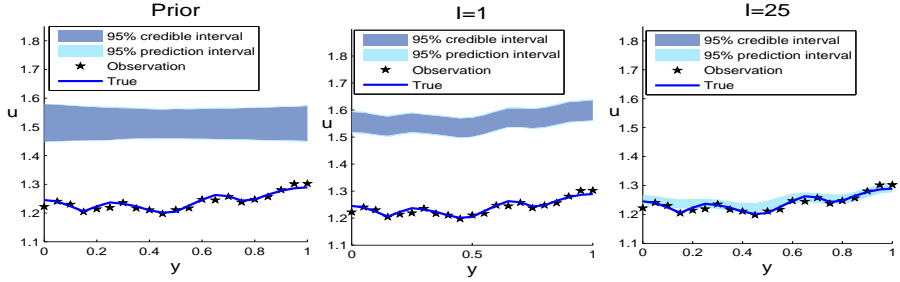}
  \caption{95\% prediction intervals, 95\% credible intervals, observations and truth values for $u((0.6,y);1)$ in different iterations.}\label{C_prey}
\end{figure}

To access the prediction using the posterior ensemble, we compute the $95\%$ credible and prediction intervals for the models response at $u((0.6,y);1)$, as shown in Figure \ref{C_prey}. We note that the true values and observations lie outside of both the credible and prediction intervals for the prior ensemble. The prediction and credible intervals are close to the true values and observations via the iterations.
The credible interval and prediction interval  become narrower as more iterations are implemented. This implies  that  the uncertainty of the parameter  $\bt$ decreases against the iterations, and the difference  between the model fit and predictions decreases. The observations data are almost concentrated in the prediction interval in the final interation.

%%%%%%%%%%%%%%%%%%%%%%%%%%%%%%%%%%%%%%%%%%%%%%%%%%%%%%%%%%%%%%%%%%
\subsection{Identify fracture structure and fractional derivative order}\label{fracture}
In this subsection, we consider the time fractional diffusion model 
\[
^{c}D_t^{\alpha} {u}-\nabla\cdot(a(x)\nabla u)=f(x,t),\quad \quad x\in\Omega, \quad t\in (0, T],
\]
where the boundary conditions are mixed and the same as in Subsection \ref{channel}.  Here we consider  the Caputo fractional derivative, i.e.,
\[
  ^{c}D_t^{\alpha} {u}=\frac{1}{\Gamma(1-\alpha)}\int_0^{t}(t-\tau)^{-\alpha}\frac{\partial u(x, \tau)}{\partial \tau}d\tau,\qquad 0<\alpha<1 ,
\]
where $\Gamma(\cdot)$ is the Gamma function. The truth  permeability field $a(x)$ is given by
\[
  a(x)=
   \begin{cases}
    10000,& x=0.3,\quad0.4\leq y\leq0.8,\\
     1, & \text{otherwise}.
  \end{cases}
\]

For this example, both the structure of the fracture and the fractional derivative $\alpha$  are unknown. We only have the prior information that the fracture parallels to $y$ axis of spatial domain. The location and the length of the fracture are unknown. The fracture in a two-dimensional permeability field is reduced by a one-dimensional line segment. Thus we describe the fracture by the midpoint coordinate $(x_0, y_0)$ and the  length $L_0$. Thus we need to identify the unknown parameter ${\bt}=(\alpha,x_0,y_0,L_0)^T$. We have the prior information $\bt\in[0,1]^4$,  which is assume to be  the  uniform distribution in $[0,1]^4$. The parameter  values are uncontrollable during the iteration of IES, the samples may run out of the  interval. Thus, a transformation is  used to overcome this problem. To this end, we make a bijective map $\bb O$: $\bb R^4\rightarrow [0,1]^4$, i.e.,
\[
{\bt}^j=\frac{1}{2}+\frac{1}{\pi}\arctan({\bm q}^j),\quad j=1,\cdots,N_e.
\]
At the  $l$-th iteration step in IES,  we use 
\[
 \left\{
 \begin{aligned}
 &\bt_l={\bb O}({\bm q}_l) \\
  &{\bm y}_l=g({\bt}_l)+\bm\varepsilon.
  \end{aligned}
\right.
\]

\begin{figure}
  \centering
  \includegraphics[width=5in,height=1.3in]{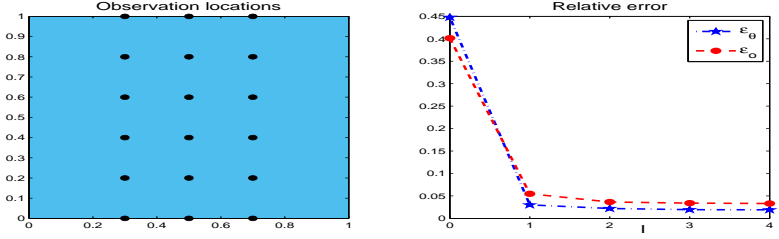}
  \caption{The observation locations (left), the relative error of $\bt$ (blue dash line) and the end points (red dash line) for the fracture via the iterations (right).}\label{frc}
\end{figure}
The truth midpoint $(x_0,y_0)=(0.3,0.6)$ and truth  length $L_0=0.4$. Then the coordinate of the true end points is ${\bm o}^{tr}=(0.3,0.4,0.3,0.8)^T$. For the source term, we take $f=10$. The end time $T=5$.  Let truth ${\bt}^{tr}=(0.7,0.3,0.6,0.4)^T$. The forward model is solved  on a uniform $100\times100$ grid with $\Delta t=0.05$ and observations are obtained by solving the same grid with time step $\Delta t=0.1$.   DFM model \cite{Efendiev2015hierarchical} is used to treat  the fracture model. The observation data are taken at time instance $t=5$, where the observation locations are distributed on the uniform $3\times 6$ grid of the domain $[0.3, 0.7]\times [0, 1]$ as shown in Figure \ref{frc} (left). The standard deviation of noise $\sigma=0.03$ and the scale $\rho=10$ in equation(\ref{weight}). We use GMM to model the prior and posterior.  The initial $\{\bm\mu_i,\bm\Sigma_i\}_{i=1}^k$ is arbitrarily given and $\pi_i=\frac{1}{k}$. The initial $k$ is set as $k_{max}$. We take $k_{min}=2$ and $k_{max}=5$.

Let $\bm o$ denote the vector, which consists of the two end points of the fracture. To measure the discrepancy between the truth and the estimate parameter, the relative errors of $\bt$ and $\bm o$ are defined by
\begin{equation*}
\varepsilon_{\theta}:= \frac{\|\widehat{\bt}-{\bt^{tr}}\|}{\|{\bt^{tr}}\|}, \quad \varepsilon_o:=\frac{\|\widehat{\bm o}-{\bm o^{tr}}\|}{\|{\bm o^{tr}}\|}, 
\end{equation*}
where $\widehat{\bt}$ is given by equation (\ref{unknown_me}).  These relative errors are shown in Figure \ref{frc} (right), where we  find that  the relative errors are large at the beginning and decrease significantly as
the observation information is incorporated into IES.  This figure also shows the convergence of  GMM-based IES-IS algorithm.

\begin{figure}
  \centering
  \includegraphics[width=5in,height=1.4in]{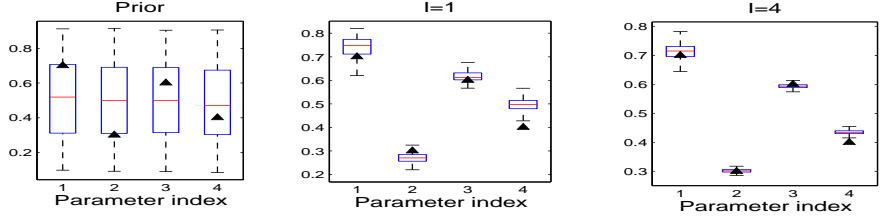}
  \caption{Medians (red solid line), reference (black filled triangle), the edges of the box (blue unfilled rectangle) are the 25th (bottom) and 75th (top) percentiles, and $95\%$ credible intervals (black dash line) for $\bt$ via the iteration, respectively.}\label{fracb}
\end{figure}

Figure \ref{fracb} depicts the 95\% uncertainty bands of parameter $\bt$ using the posterior ensemble samplers at different iterations.
We see the medians of the prior ensemble are  in the middle of the blue rectangle and the credible intervals are large. This is due to that the samplers are drawn from the Gaussian distribution at the beginning. The credible intervals gradually become narrow as more iterations are used. This implies that the uncertainty of the unknowns decreases in the IES process. In the last iteration, the 95\% credible intervals are tight and the reference values are included in the credible intervals except for the length of fracture $L_0$. This also shows  the support of posterior is only a small portion of the support of prior ensemble.

\begin{figure}
  \centering
  \includegraphics[width=5in,height=1.4in]{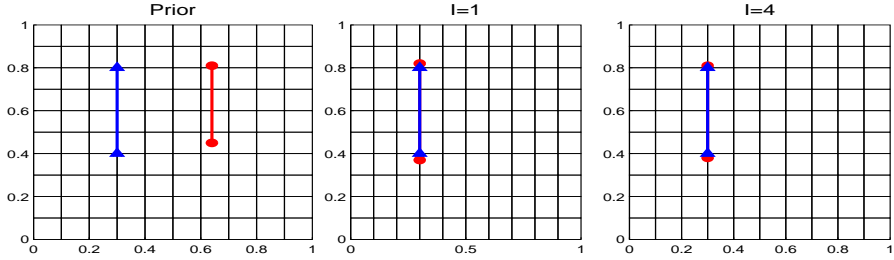}
  \caption{The inversion  of the fracture via the iteration steps. The red segment  is the estimate fracture and the blue segment  is the truth fracture.}\label{frac0}
\end{figure}

Figure \ref{frac0} depicts the recovered fracture against the iteration procedure.
By the figure,  the discrepancy between    the prior  and truth fracture is very large, which may lead big difference for the simulated  state. When the information of the observations enters the inference by IES-IS, the estimate fracture gets  close to the truth. Compared with the truth fracture, we  find the estimate of the middle point is more accurate than the  length. This is due to that we use DFM model, where the fracture must match with the mesh grid of the forward model. As expectation, the uncertainty of the fracture decreases as the iteration moves on. 

\begin{figure}
  \centering
  \includegraphics[width=5in,height=1.4in]{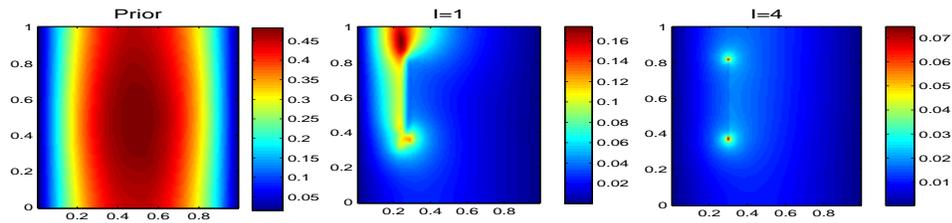}
  \caption{The posterior standard deviations of the simulated state $u((x,y);5)$ at different iteration steps.}\label{frac_err}
\end{figure}

We also plot the posterior standard deviations of the simulated state $u((x,y);5)$ at different iteration steps in Figure \ref{frac_err}. We see  that the standard deviation is small around the boundaries corresponding to  the deterministic Dirichlet boundary condition. For the prior ensemble, the standard deviation is large  in the most part of  domain. Then  the standard deviation decreases as uncertainty is reduced with respect to IES iterations. In the inversion process, the region of the large standard deviation becomes narrow and is mostly  concentrated around the fracture.  In  the final posterior ensemble, the standard deviation in the whole physical domain is very small and the uncertainty  is almost only around the  endpoints of the fracture. 

%%%%%%%%%%%%%%%%%%%%%%%%%%%%%%%%%%%%%%%%%%%%%%%%%%%%%%%%%%%%%%%%%%%%%%%%%%%%%%%%%%
\section{Conclusion}
We proposed an ensemble-based implicit sampling for handling non-Gaussian priors.
 In the approach ,  iterative ensemble smoother (IES)  has been coupled with implicit sampling (IS).  IES can efficiently  provide  an approximation to the MAP point of the posterior and the inverse of Hessian matrix.
 IES avoids the explicit computation of Jacobian matrix and Hessian matrix for the optimization problem in the underlying inverse problem.  This computation is usually challengeable in high dimension spaces.
 Then an implicit map was constructed by the  the MAP point and the  Cholesky factorization of the Hessian matrix.
 IS was used to identify  a high probability region and obtain the importance samples.
 The proposed method was extended to non-Gaussian priors where DCT and GMM are used to characterize prior. This significantly improved the applicability of the conventional implicit sampling.
 
 We applied the proposed sampling method to inverse problems of subsurface flows and anomalous diffusion models in heterogeneous porous media.  
 The  ensemble-based implicit sampling was used to effectively recover channel structures and fracture structures in porous media for these models.

\bibliographystyle{siamplain}
\bibliography{DC-En}

\end{document}